\documentclass[12pt,reqno]{amsart}
\usepackage{tikz,dsfont}
\usepackage{amsmath, amsthm,amssymb, graphicx, mathrsfs, hyperref}
\usepackage{pdfpages}
\usepackage{verbatim} 
\usepackage{enumitem}

\usepackage{hyperref}

\usepackage{xcolor}
\usepackage{soul}

 \usepackage{float}
\newtheoremstyle{dotless}{}{}{\itshape}{}{\bfseries}{}{ }{} 
\theoremstyle{dotless}

\newtheorem{theorem}{Theorem}[section]

\newtheorem{lemma}{Lemma}

\newtheorem{prop}{Proposition}
\newtheorem{conj}{Conjecture}

\newtheorem{corollary}{Corollary}

\usepackage{color}
\definecolor{red}{rgb}{1,0,0}
\definecolor{blue}{rgb}{.2,.6,.75}
\definecolor{green}{rgb}{.4,.7,.4}
\newcommand{\acom}[1]{{\color{red}{ #1}} }

\title[A sum formula and abc]{A  sum formula for finite subspaces in function fields and extensions of the $abc$-conjecture}
 
\author{Andrew Granville}

\begin{document}

\begin{abstract}
Let $K$ be a finite extension of the rational function field $k(t)$. 
We show that there is an appropriate extension of the definition of each valuation $v$ of $K$ such that
for any finite dimensional subspace $T$ of $K$ (considered as a vector space over $k$) of dimension $n$, we have the sum formula
 \[
 \sum_v v(T) = \binom{n}2 (2g-2)
\]
where $g$ is the genus of $K$. Then for any finite subset  $X$ of $T$ we define $v(X)$ appropriately so that $v(X)\leq v(T)$ and ask when the inequality $ \sum_v v(X) \leq \binom{n}2 (2g-2)$ is an  equality, particularly for sets $X=\{ x_0,\dots,x_n\}$ where $x_0+\cdots +x_n=0$. The $n=2$ case corresponds to classifying Bely\u{\i} maps, a well-known difficult problem that  is still far from fully resolved, and here we begin investigations into the $n>2$ case. This gives us some insight into extending the $abc$-conjecture to arbitrarily many variables.
\end{abstract}

\maketitle

\section{Introduction and results}

Let $k$ be an algebraically closed field of characteristic zero and  $K$ be a finite extension of the rational function field $k(t)$; the  valuations on $K$  correspond to points on a nonsingular curve $C$ defined over $k$ (the valuations  are 0 on $k^*$ but are not identically 0, and their value group is isomorphic to $\mathbb Z$ as an ordered group). The  valuation $v$ at a point is given by the order of vanishing of the Laurent
expansion in a local parameter $z_v$, and so valuations  are non-archimedean, additive and integer-valued
(valuations with $v(z_v)>0$ are ``finite'', while those with $v(z_v)<0$ are ``infinite'' and there are only finitely many infinite valuations).   
Finite valuations of a function field $K$ correspond to prime ideals (or places) of its integral closure.
The (well-known) \emph{sum formula} is given by 
\begin{equation} \label{sum.formula}
\sum_v  v(f)=0 \text{ for all non-zero } f \in K,
\end{equation}
where the sum is over all valuations $v$ of $K$.  For any  subset $X=\{ x_1,\dots,x_n\}\subset K^*$ with $n\geq 2$ we define the \emph{projective height} as
\[
H(X):=- \sum_v \min\{ v(x): x\in X\}
\]
which is $\geq 0$, and $H(X)=0$ if and only if   $x/y$ is a unit for all $x,y\in X$. Now if $x$ is non-zero then $v(x)\ne 0$ for only finitely many valuations $v$ and so determining $H(X)$ when $X$ is finite involves only finitely many valuations.  If $v(x)=0$ for all $v\not\in S$, a set of valuations, then we call $x$ an \emph{$S$-unit}. Typically, in applications, $S$ is taken to contain all the infinite valuations.

The sum formula is a useful tool in algebraic number theory, often linked to investigations of heights in Diophantine equations, and various height inequalities are conjectured in the literature (like the $abc$-conjecture, and Vojta's conjectures), both in the function field setting and in number fields.

 Suppose that $x_1,\dots,x_n$ is a basis for  the $n$-dimensional subspace $T$ of $K$, where $T$ is to be thought of as a vector space over $k$.
We will  extend the definition of  $v(\cdot)$ to any such subspace $T$   in the next subsection
so that $\{ v(T)\}_v$ reveals the ramification structure of the map 
\[
\Phi(t)=(x_1 :\dots:x_n)\in \mathbb P^{n-1}(K) \text{ with } \Phi: k\to \mathbb P^{n-1}(k). 
\]
We obtain the following
\emph{sum formula for finite-dimensional subspaces}: 

\begin{theorem}  \label{thm: main} 
If $K$ is a function field of genus $g$ over $k$, and $T$ is an $n$-dimensional subspace of $K$
(as a vector space over $k$) then
\begin{equation} \label{sum.formula 4 S}
 \sum_v v(T) = \binom{n}2 (2g-2)
\end{equation}
where the sum is over all valuations $v$ of $K$.  
 \end{theorem}

There are formulas in the literature which are similar but formulated differently (typically in terms of Weierstrass points on the nonsingular curve $C$ mentioned above),  such as \cite[(3)]{Vo2}, the  ``Brill-Segre formula'' in Voloch \cite{Vo},
or the ``Pl\"ucker formula'' in Eisenbud and Harris \cite{EH}, who in turn reference Weyl and Piene.
The dimension $2$ case is equivalent to the Riemann-Hurwitz formula applied to rational maps $\phi: \mathbb P^{1}(k)\to \mathbb P^{1}(k)$, where $T$ is generated by the numerator and denominator of $\phi$ (as in \cite{Si}).

Calculating all the $v(T)$-values may not be practical in many applications. Often we are presented with
a finite subset $X\subset T$ which spans $T$ as a vector space and we want bounds in terms of the set $X$ rather than $T$. 
 For any finite subset $X\subset T$ we will define
$v(X)$ below, easy to calculate in terms of the (finite number of) elements of $X$, where $v(X)\leq v(T)$ for all valuations $v$, and so we can deduce the following:

\begin{corollary} \label{corr: sum of the x_i}
Let $X \subset K^*$, a function field of genus $g$ over $k$, where $X$ generates an $n$-dimensional subspace $T$.
Then
\begin{equation} \label{sum.formula 4 x}
 \sum_v v(X)\leq\binom{n}2 (2g-2)
\end{equation}
where the sum is over all valuations $v$ of $K$.   We can restrict the summation here to ``$v\in S$'' provided every element of $X$ is an $S$-unit.  %
  \end{corollary} 

 The main novelty in this article, as compared with most of the existing literature, is \eqref{sum.formula 4 x} because it develops the theory for arbitrary finite spanning subsets  $X$ of $T$ rather than for solutions to  equations like the Diophantine equation  
\begin{equation} \label{sum x's}
x_0+x_1+\cdots+x_n=0
\end{equation}
where  the $x_i$ generate an $n$-dimensional subspace $T$,  and we take $X=\{ x_0,\dots,x_n\}$. We also obtain new results in this case.

Given any integers $d_0,d_1,\dots, d_n$, reorganize them so that $d_{\sigma_0}\leq d_{\sigma_1}\leq  \cdots \leq d_{\sigma_n}$ and define
\[
\mathcal E(d_0,\dots,d_n)=\mathcal E(\{d_0,\dots,d_n\}) := \sum_{j=1}^{n} \bigg( d_{\sigma_j}\  - \max_{1\leq i\leq j} (d_{\sigma_i}-(i-1) ) \bigg) 
\]
which, we will show, is in $\in [0,\binom n2]$, and equals $0$ if and only if the $d_i$ are all equal.   

\begin{theorem} \label{corr: sum of the x_i=0}
Suppose that we have a solution to \eqref{sum x's}
where  the $x_i$ generate an $n$-dimensional subspace $T$,  and we take $X=\{ x_0,\dots,x_n\}$.  Then
\begin{equation} \label{eq: ht ineq2}
H(X)   \leq   \sum_{v\in S}  \mathcal E_v(X) + \binom{n}2 (2g-2) 
\end{equation}
where the elements of $X$ are all $S$-units, and we define $\mathcal E_v(X):= \mathcal E (v(x_0),\dots,v(x_n))$.
\end{theorem} 

For a solution to  \eqref{sum x's}, the two inequalities \eqref{sum.formula 4 x} and \eqref{eq: ht ineq2} imply each other, but this latter form may be easier for applications as it is completely explicit.

Since each $ \mathcal E_v(X) \leq \binom{n}2$ we deduce 
\begin{equation} \label{sum.formula 4 H}
 H(X)=H(T)\leq\binom{n}2 (2g-2+\# S);
\end{equation}
from \eqref{eq: ht ineq2}.
The inequality \eqref{sum.formula 4 H} is given in \cite[Corollary I]{BM} and \cite[Theorem 4]{Vo2} where it is claimed to be originally due to Mason (though he only gave the case  $n=2$ in \cite{Ma}), and also as \cite[Corollary]{Za}. The sharper inequality \eqref{eq: ht ineq2} appears to be new and is probably as much as can be squeezed out of our methods without even more complicated definitions.  This is most intuitively applied to Diophantine equations involving polynomials:

\begin{corollary} \label{corr: sum of the x_i=0. C[t]}
Suppose that we have a solution to \eqref{sum x's}
where  the $x_i\in \mathbb C[t]$ have no common root and generate an $n$-dimensional subspace $T$ (so no subsum can be $0$),  and we take $X=\{ x_0,\dots,x_n\}$. 
Then
\begin{equation} \label{eq: Sum Bound}
  \max_i \deg x_i  \leq  \sum_{v:\ v(x_i)\ \text{ not all equal }}     \mathcal E_v(X)  -2\binom{n}2 .
\end{equation}
Here $\min_i v(x_i) =0$ if $v=v_{t-\alpha}$, and $\min_i v(x_i) =-\max_i \deg x_i $ if $v=v_\infty$.
\end{corollary} 

For a solution to  \eqref{sum x's}, the  inequality \eqref{sum.formula 4 x} in $\mathbb C[t]$ and \eqref{eq: Sum Bound} imply each other; it is important to note that the sum over $v$ in each inequality includes the infinite valuation.
We will study equality in these equations in section  \ref{sec: OPP}.

If the $x_i$ are pairwise coprime then each $\mathcal E_v(X)= \min\{ \max_i v(x_i),n-1\}\leq n-1$ for $v\ne v_\infty$, and 
$\mathcal E_{v_\infty}(X)\leq \binom{n}2$. Therefore
\begin{equation} \label{eq: Coprime}
 \max_i \deg x_i  \leq  (n-1) \#\{\alpha\in \mathbb C:\ x_i(\alpha)=0 \text{ for some } i\}       -\binom{n}2 .
\end{equation}
We will see applications of this to Diophantine equations in section \ref{sec: Applications} below.

In the $n=2$ case, $\mathcal E_v(X) =1$ whenever $v(x_0), v(x_1), v(x_2)$ are not all equal and $\mathcal E_v(X) =0$ otherwise, so we obtain the usual \emph{$abc$-theorem for polynomials}, which in turn suggests the $abc$-conjecture. It has been a challenge to guess at the right exponents on primes in the $abcd$-conjecture and beyond.  We feel more ready to guess now, based on an analogy to Corollary \ref{corr: sum of the x_i=0. C[t]}:
\medskip

\noindent \textbf{The extended $abc$-conjecture} \emph{ Suppose that we have a solution to \eqref{sum x's} in integers $x_0,\dots,x_n$ with no common factor, and that no subsum vanishes. Then
\[
\max_i |x_i| \ll_\epsilon \bigg(  \prod_{p|x_0\cdots x_n}  p^{\mathcal E_p(x_0,\dots,x_n)}\bigg)^{1+\epsilon}
\]
where $\mathcal E_p(x_0,\dots,x_n)=\mathcal E( v_p(x_0),\dots, v_p(x_n) )$, and $v_p(x)$ denotes the exact power of $p$ dividing the integer $x$.  
}
\medskip

The results in section \ref{sec: OPP} suggest that  if true the extended $abc$-conjecture is best possible in many circumstances.    The formulation of the extended $abc$-conjecture is complicated  because of the need to evaluate the $\mathcal E_p$; it is desirable to work with a more easily understood conjecture.
Define the \emph{radical}, $r(x)=\prod_{p|x} p$.
Now if each $d_i\geq 0$ then $\mathcal E(d_0,\dots,d_n)\leq (n-1)\#\{ i: d_i\geq 1\}$ and so
for each prime $p$ we have
\[
\mathcal E_p(x_0,\dots,x_n) \leq \ (n-1) \# \{ i: p|x_i\}.
\]
Therefore the extended $abc$-conjecture implies that 
\[
\max_i |x_i|  \ll_\epsilon \bigg(  \prod_{i=0}^n  r(x_i) \bigg)^{n-1+\epsilon}.
\]
If one has an example close to equality  here only then the $x_i$ are (more-or-less) pairwise coprime
and most of the primes  divisors divide their $x_j$ to a power $\geq n-1$.
Since each $\mathcal E_p(x_0,\dots,x_n)\leq \binom n2$, the extended $abc$-conjecture also implies  that
\[
\max_i |x_i|  \ll_\epsilon  r(x_0\cdots x_n)^{\binom n2+\epsilon};
\]
moreover if one has an example close to equality  here then most of the prime divisors  divide  all but two of the $x_j$, and then to all different (positive) powers.
Browkin and Brzezi\'nski \cite{BB} conjectured the sharper upper bound 
\[
\max_i |x_i|  \ll_\epsilon r(x_0\cdots x_n)^{2n-3+\epsilon}
\]
and showed that   the exponent $2n-3$ cannot be reduced. Notice that $2n-3=\binom n2$ for $n=2$ and $3$, so which conjectured exponent is  correct can perhaps only be determined with extreme examples for $n=4$.

In Proposition \ref{prop BB} we will show that Browkin and Brzezi\'nski's conjecture implies that 
for any solution to  \eqref{sum x's} where  the $x_i\in \mathbb Z[t]$ have no common root and generate an $n$-dimensional subspace $T$ we have
\[
 \max_i \deg x_i  \leq  (2n-3)  (\# \{\alpha\in \mathbb C:\  x_i(\alpha)=0 \ \text{ for some } i\} -2) .
\]
If true this is best possible as we can exhibit an example achieving equality but we have no compelling evidence that this is true.
If this inequality is false then Browkin and Brzezi\'nski's conjecture is also false.

In \cite[section 8]{Gr} we suggested that 
\[
\max_i |x_i|  \ll_\epsilon r(x_0\cdots x_n)^{1+\epsilon}
\]
outside of finitely many subvarieties (which was again conjectured by Looper \cite{Lo} with a nice application). 
It is difficult to relate this to polynomial identities, because in section \ref{sec:  Inf Id} we will construct 
infinitely many identities with implicit exponent $<\epsilon$, unrelated by M\"obius transformations, and yet all on the same subvariety.

It seems to be of particular interest to determine whether there are examples of solutions to \eqref{sum x's} in 
$\mathbb C[t]$  for which we get equality in \eqref{sum.formula 4 x}, where $X=\{ x_0,\dots,x_n\}$.  
This is equivalent to determining those maps  
\[
\Phi(t)=(x_1(t):\dots:x_n(t))\in \mathbb P^{n-1}(\mathbb C(t))  \text{ with } \Phi: \mathbb C\to \mathbb P^{n-1}(\mathbb C)
\]
for which all the ``ramification structure'' lies above
\[
(1,0,\dots,0), \dots, (0,\dots,0,1,0,\dots,0), \dots, (0,\dots,0,1) \text{ and } (1,1,\dots,1).
\]
which correspond to the roots of $x_1(t),\dots ,x_n(t)$ and $x_0(t)$, respectively.
In the $n=2$ case this ``ramification structure''  corresponds to classifying Bely\u{\i} maps, a well-known difficult problem that can be approached in many different and interesting ways (see \cite{SV} and section \ref{sec: Belyi}), but is still far from fully resolved.
For $n>2$ the ``ramification structure'' is more complicated and involves the $\{ x_j(t)\}_{0\leq j\leq n}$ being divisible by different specific powers of $(t-\alpha)$ for those $\alpha$ with $v_{t-\alpha}(T)>0$; there seems to be very little relevant literature.  We will present infinitely many examples of equality for $n=3$ in section \ref{sec: uplift of n=2} and for each $n$ in sections \ref{sec: dth power} and \ref{sec:  power + roots}, but defer a broader study to the further paper \cite{AI}.

 Explorations via AI are starting to have an impact on research mathematics. Researchers at Axiom Math worked with their AxiomProver  to help develop an earlier version of this preprint in several ways. 
 First they used AxiomProver to verify the proof of Theorem \ref{thm: sumformula for C[t]} in Lean; AxiomProver automated a translation of my ``natural-language'' proof into machine-verifiable code.
 In the earlier version of the paper we had not established Theorem \ref{corr: sum of the x_i=0} and Corollary \ref{corr: sum of the x_i=0. C[t]} but rather made do with weaker inequalities like \eqref{sum.formula 4 H},
  \eqref{eq: Zannier.1c} and \eqref{eq: Zannier.7} (in what was then section 6); AxiomProver  autonomously
sharpened these inequalities obtaining   Corollary \ref{corr: sum of the x_i=0. C[t]}  with $n=3$, where
 $\mathcal E(0,d_2,d_3)$ was evaluated as in section \ref{sec: E eval}. It also looked further into determining examples of  the polynomial identity suggested in section \ref{sec: E eval}.  A discussion of the contributions of AxiomProver can be found in  Appendix \ref{sec:AxiomProver}.
 
We will first prove our new sum formula \eqref{sum.formula 4 S} in $\mathbb C[t]$, the simplest case, and then modify the proof for general function fields. In the next section we explain our extensions of valuations, and prove various straightforward properties.

\section{Extending  valuations} \label{sec: Extending valuations}

\subsection{Extending  valuations to finite dimensional subspaces} \label{sec: Extending valuations 2 subspaces}
 Fix a valuation $v$ of $K$. We will  show that there exists a basis $y_{v,1},\dots, y_{v,n},$ of $T$ for which
\begin{equation} \label{growing ys}
 v(y_{v,1})<v(y_{v,2})<\dots <v(y_{v,n}).
 \end{equation}
 We call this an \emph{ordered basis for $T$ with respect to $v$}.
 Let $T_j:=  \{ y\in T: v(y) \geq v(y_{v,j})\}$ which is a subspace of $T$ of dimension $n+1-j$ with $T_n\subsetneq T_{n-1}\subsetneq \cdots \subsetneq T_1=T$.

 We define
 \[
v(T) := \sum_{j=1}^n v(y_{v,j}) - \binom{n}{2}. 
\]
If $y\in T$ then $v(y)=v(y_{v,i})$ for some $i$ which implies that $v(T)$     is independent of the choice of basis (see Lemma \ref{ChangeOfBasis2}).  We can therefore define 
\[
H(T):=- \sum_v \min\{ v(y): y\in T\}
\]
and so if  $X\subset T$ then $\min\{ v(x): x\in X\} \geq  \min\{ v(y): y\in T\}$ and thus
\[
H(X) \leq H(T) .
\]
We obtain equality here if $x_1,\dots ,x_n$ generates $T$ since then $\min\{ v(x): x\in X\} = \min\{ v(y): y\in T\}$ for each $v$.

 Now \eqref{growing ys} implies that $v(y_{v,j})\geq v(y_{v,1})+j-1$ for each $j$ and so $v(T)\geq nv(y_{v,1})=n\min_{y\in T} v(y)$.    
 If $v(y_{v,j})> v(y_{v,1})+j-1$ for some $j$ (which is equivalent to $v(T)>n\min_{y\in T} v(y)$) then $v$ is a \emph{Weierstrass point} (see \cite{Vo2}).

There are only finitely many finite valuations $v$ for which $v(y)<0$  for some $y\in T^*$, for if 
$y_1,\dots,y_n$ is a basis for $T$ then some $v(y_j)<0$ and so $v$ corresponds to a prime ideal dividing the denominator of $y_j$; and there can be only finitely many for each $j$.  Therefore $v(T)\geq 0$ aside from some finite set of valuations (which may include all the infinite valuations), and so, since each $v(T)\in \mathbb Z$ we deduce from  \eqref{sum.formula 4 S} then $v(T)=0$ for all but finitely many $v$.  Hence there are only finitely many 
Weierstrass points.
   
 If the elements of $T$ are integral with respect to all finite valuations $v$
 then,  for every infinite valuation $v$, $v(y)\leq 0$ for all $y\in T$ and so $v(T)\leq -2 \binom{n}{2}$.

\subsection{Integral elements} \label{sec: Integral}
We can assume, without loss of generality, that $T$ is a subspace of the ring of integers of $K$ (so that the elements of $T$ are integral with respect to all finite valuations $v$). To see this, determine any basis $y_1,\dots,y_n$ for $T$ as an $n$-dimensional subspace of $K$ over $k$.  Select $h\in K$ so that $hy_j$ is in the ring of integers of $K$ for each $j$ and then $h\cdot T:=\{ hy: y\in T\}$ is a subspace of the ring of integers of $K$ (which  forms a finite dimensional vector space over the ring of integers of $k$).
Then $v(hT) = n v(h) +v(T),\ v(hX) = n v(h) +v(X)$ and $\min_{z\in hT} v(z) = v(h) + \min_{y\in T} v(y) $ so that 
$\sum_v v(hT)= \sum_v v(T),\ \sum_v v(hX)= \sum_v v(X)$ and $H(hT)=H(T)$ using the sum formula. 
This implies that 
we can deduce Theorem \ref{thm: main} and Corollary \ref{corr: sum of the x_i} from the same result restricted to subspaces of the ring of integers, using \eqref{sum.formula} with $f=h$.  

So now assume that $T$ is a subspace of the ring of integers of $K$. Let $G$ be the greatest common divisor of the elements of $T$, which is an ideal.
If $K$ is a PID then we may divide $T$ through by $g$ where $G=(g)$ so that  $v(y_{v,1})=0$ for each finite place $v$.
 In that case $v(T)\geq 0$, and $v(T)>0$ if and only if $v(y)\geq n$ for some $y\in T$.  In non-PIDs we may select a prime ideal $P$ in the same ideal class as $G$, so that $G\overline{P}=(h)$ and $P\overline{P}=(p)$ are principal ideals, so scaling $T$ by the factor $p/h$ leads to a ring of integers, with greatest common divisor $P$, a prime ideal. Here $v(y_{v,1})=0$ for each finite place $v$, except the place $v_P$ corresponding to $P$, whence $v_P(y_{P,1})=1$.

\subsection{Extending  valuations to subsets of a finite dimensional subspace} \label{sec: Extending valuations 2 sums}

 Suppose that $X \subset K$  spans the $n$-dimensional subspace $T$. Let 
 \[
 X_{v,j}=\{ x\in X: v(x)\geq v(y_{v,j})\}\subset T_j
 \]
for each $j$, so  that $X_{v,1}=X$. For
$j=n,n-1,\dots, 1$ select a linearly independent set $Y_{v,j}\subset X_{v,j}$, with $Y_{v,j+1}\subset Y_{v,j}$ for $j\leq n-1$, such that  $\langle Y_{v,j}\rangle_k = \langle X_{v,j}\rangle_k$
(that is, $Y_{v,j}$ generates the same subspace of $T$ as the elements of $X_{v,j}$). We label the $n$ elements of $Y_{v,1}$ so that
\[ v(x_{v,1})\leq \cdots\leq v(x_{v,n}).\]
 Since  $\langle Y_{v,1}\rangle_k = \langle X_{v,1}\rangle_k= \langle X \rangle_k=T$ we have $v(x_{v,1})=v(y_{v,1})$.
This implies that 
\[
H(X)=- \sum_v \min_{x\in X} v(x)= - \sum_v v(x_{v,1}) = - \sum_v v(y_{v,1}) =H(T).
\]

By construction we see that for any choice of $x_1,\dots, x_n\in X$ ordered so that $v(x_1)\leq \cdots \leq v(x_n)$ we have 
$v(x_j)\leq v(x_{v,j})$ for each $j$. Therefore if $U\subset X$ and $U$ generates $T$ also then 
$v(u_{v,j})\leq v(x_{v,j})$ for each $j$.

Now  
\begin{equation} \label{eq: pure vS2}
v(T) = \sum_{j=1}^{n} v(x_{v,j}) + \Delta_v(X)  -\binom n2
\end{equation}
where the \emph{deficit} is defined by
 \[
\Delta_v(X) =\Delta_v(x_{v,1},\dots,x_{v,n}):= \sum_{j=1}^{n} v(y_{v,j}) - \sum_{j=1}^{n} v(x_{v,j}) \quad (\geq 0).
\]

Define, for  any sequence of integers  $d_1\leq\cdots \leq d_n$, the \emph{minimum deficit function}
\[
\Delta(d_1,\dots,d_n):=\sum_{i=1}^n (D_i-d_i)
\]
where  $D_1<D_2<\cdots<D_k$ is a sequence of integers selected minimally so that $D_j\geq d_j$ for each $j$.  

\begin{lemma} \label{lem: Deficit} Let $d_0\leq d_1\leq\cdots \leq d_n$ be a given sequence of integers. Select $D_1<D_2<\cdots<D_k$ minimally so that $D_j\geq d_j$ for each $j$.  Then 
\[
D_j= \max_{i\leq j} (d_i+j-i ) \text{ and } 0\leq \Delta(d_1,\dots,d_n) \leq \binom n2,
\]
so that $\mathcal E(d_0,\dots,d_n)  = \binom n2-\Delta(d_1,\dots,d_n)  \in [0,\binom n2]$.
\end{lemma}

We also note that $\Delta(d_1,\dots,d_n) =\binom n2$ if and only if $d_1=\cdots=d_n$,
and $\Delta(d_1,\dots,d_n) =0$ if and only if $d_{j-i}\leq d_j-i$ for all $i\leq j$.

\begin{proof}
Now each $D_j \geq D_{j-1}+1$ and so if $j>i$ we have $D_j\geq D_i+(j-i)\geq d_i+j-i $, from which we can deduce that 
$D_j\geq\max_{i\leq j} (d_i+j-i )$, and this sequence satisfies the hypothesis, which implies the first result. 

Now  $d_j\leq D_j\leq \max_{i\leq j} d_i +(j-1)=d_j+j-1$, and so
$0\leq \Delta(d_1,\dots,d_n) \leq \binom n2$. 

With the $d_i$'s ordered we have, by definition,
\begin{align*}
\mathcal E(d_0,\dots,d_n) &= \sum_{j=1}^{n} \bigg( d_j - \max_{1\leq i\leq j} (d_i-(i-1) ) \bigg) \\
&= \sum_{j=1}^{n} ( d_j +j-1- D_j  ) = \binom n2-\Delta(d_1,\dots,d_n) .\qedhere
\end{align*}
\end{proof}

Define 
\[
v(X):= \sum_{j=1}^{n} v(x_{v,j}) + \Delta( v(x_{v,1}),\dots,v(x_{v,n})) -\binom n2.
\] 
We observe that if $v(x)=0$ for all $x\in X$ then each term in the sum is $0$ and sum of the last two terms equals $0$, so that $v(X)=0$. Hence $v(X)=0$ for all but finitely many valuations $v$.

\begin{lemma}\label{lem: size of v(x)} 
If $U\subset X$ then $v(X)\geq v(U)$ for every valuation $v$; in particular $v(T)\geq v(X)$.     
We also have $v(X)\geq nv(y_{v,1})$ with equality if and only if $0\leq v(x_{v,i})-v(y_{v,1})\leq i-1$ for each $i$.
\end{lemma}

\begin{proof}
Now if  $d_j=v(x_{v,j})$ for each $j$ then $v(y_{v,j})\geq  D_j=\max_{i\leq j} (v(x_{v,i})+j-i )$ for each $j$ (by the minimality of the $D_j$) and so 
\[
\Delta_v(X)\geq  \Delta( v(x_{v,1}),\dots,v(x_{v,n})) .
\]
Therefore $v(T)\geq v(X)$ by their definitions. We also deduce that 
\[
v(X)= \sum_{j=1}^{n} \max_{i\leq j} (v(x_{v,i})+j-i )-\binom n2,
\]
so if $U\subset X$ then    
\[
v(X)= \sum_{j=1}^{n} \max_{i\leq j} (v(x_{v,i})+j-i )-\binom n2\geq \sum_{j=1}^{n} \max_{i\leq j} (v(u_{v,i})+j-i )-\binom n2=v(U).
\]
Since $v(x_{v,i})+j-i \geq v(x_{v,1})+j-1$ we also obtain 
\[
v(X) \geq \sum_{j=1}^{n} (v(x_{v,1})+j-1)-\binom n2=nv(x_{v,1})=nv(y_{v,1}).
\]
  To get equality here we need 
$v(x_{v,i})+j-i \leq v(x_{v,1})+j-1$ for each $i\leq j$ and so $v(x_{v,i}) \leq v(x_{v,1})+i-1$ for each $i\leq n$, yielding the final claim.
\end{proof}

We get equality in  \eqref{sum.formula 4 x} if and only if  $\Delta_v(X) =\Delta( v(x_{v,1}),\dots,v(x_{v,n}))$ whenever $v\in S$.

\begin{proof} [Proof of Corollary \ref{corr: sum of the x_i}]   By Lemma \ref{lem: size of v(x)} 
we see that \eqref{sum.formula 4 S} implies \eqref{sum.formula 4 x}. We can restrict the sum to $v\in S$ as we saw above that 
$v(X)=0$ whenever $v(x)=0$ for all $x\in X$ (so for $v\not\in S$).    
\end{proof}

\subsection{Sums equaling zero} \label{sec: ZeroSum}
Suppose that we have a solution to \eqref{sum x's}
where  the $x_i$ generate an $n$-dimensional subspace $T$,  and we take $X=\{ x_0,\dots,x_n\}$. 

\begin{proof} [Proof of Theorem \ref{corr: sum of the x_i=0}]   
For any valuation $v$, there are at  least two $x_j$ with $v(x_j)=\min_i v(x_i)$.  Hence we may re-order the elements of $X$ so that 
\[ v(x_{v,0})= v(x_{v,1})\leq \cdots\leq v(x_{v,n}).\]
Then  for every valuation $v$ we have
\[
\sum_{j=1}^{n} v(x_{v,j}) =\sum_{i=0}^{n} v(x_{v,i}) - v(x_{v,0})=\sum_{i=0}^{n} v(x_i) - \min_j v(x_j).
\]
Now if $v\not\in S$ then  $v(x_j)=0$ for all $j\geq 1$ and so $v(x_0)\geq 0$. Therefore
\[
\sum_{v\in S} \bigg( \sum_{j=1}^{n} v(x_{v,j}) \bigg) \leq  \sum_v \bigg( \sum_{j=1}^{n} v(x_{v,j}) \bigg) = \sum_{i=0}^{n} \sum_v v(x_i) -\sum_v  \min_j v(x_j)=H(X),
\]
by the sum formula applied to each $x_i$ and the definition of height. Therefore  \eqref{sum.formula 4 x} implies that 
\begin{equation} \label{eq: ht ineq}
H(X) +  \sum_{v\in S}  \Delta( v(x_{v,1}),\dots,v(x_{v,n}))  \leq\binom{n}2 (2g-2+\# S).
\end{equation}
Subtracting $\Delta( v(x_{v,1}),\dots,v(x_{v,n})) $ from both sides for each $v\in S$ we deduce that 
\[
H(X)   \leq\binom{n}2 (2g-2)  +  \sum_{v\in S}  \mathcal E (v(x_{v,0}),\dots,v(x_{v,n}))\leq\binom{n}2 (2g-2+\# S)
\]
as each $ \mathcal E (v(x_{v,0}),\dots,v(x_{v,n}))\leq\binom{n}2$, which yields  Theorem \ref{corr: sum of the x_i=0} and \eqref{sum.formula 4 H}.
\end{proof}

\section{A sum formula for $\mathbb C[t]$ and $\mathbb C(t)$}
We begin by proving Theorem \ref{thm: main} for polynomials.
Let $T$ be a finite dimensional subspace of $\mathbb C[t]$ over $\mathbb C$ and suppose it has basis
\[
x_1(t),\dots,x_n(t).
\]
 (This is without loss of generality since we could have started with a finite dimensional subspace of $\mathbb C(t)$ over $\mathbb C$, taken a basis and multiplied through by the lowest common denominator of the basis elements to obtain a finite dimensional subspace of $\mathbb C[t]$ over $\mathbb C$, as discussed in section \ref{sec: Integral}.)

For each $\alpha\in \mathbb C\cup\{ \infty\}$ we will define $v_\alpha(x(t))$ and $v_\alpha(T)$ as follows:
For every $\alpha\in \mathbb C$   we write each $x_j(t)$ as a function of $\pi=\pi_\alpha:=t-\alpha$, so that 
\[
x_j(t)= a_j \pi_\alpha^{e_j}+\cdots \text{ (higher order terms) where } a_j\ne 0 \text{ and } e_j\geq 0,
\]
and then we define $v_\alpha(x_j(t))=e_j$.  The subspace of $\mathbb C[\pi_\alpha]$ generated by 
$x_1(t),\dots,x_{n}(t)$ over $\mathbb C$ has a basis
$y_{\alpha,1}(t),\dots,y_{\alpha,n}(t)$ with 
\[
0\leq v_\alpha(y_{\alpha,1}(t))< v_\alpha(y_{\alpha,2}(t))<\cdots < v_\alpha(y_{\alpha,n}(t))  
\]
To see this note that if $e_i=e_j$ then
$v_\alpha (a_jx_i-a_ix_j)>e_i$ so we can replace the pair $x_i,x_j$ by $x_i,a_jx_i-a_ix_j$. One keeps on doing this until the $v_\alpha(y_j(t))$ are all distinct.

We claim  that this sequence of integers is unique: Any non-zero linear element $s$ of $T$ can be written as
$s=\sum_{j=J}^{n} a_j y_{\alpha,j}(t)$ for some $J\geq 0$ with $a_J\ne 0$, and then 
$v_\alpha(s)=v_\alpha(y_{\alpha,J}(t))$. Therefore there are only these $n$ possible integers that can be values of $v(s), s\in T^*$, and so any two sets of $n$ distinct valuations of $T$ must be the same.
 
This implies that 
\[
v_\alpha(T) := \sum_{j=1}^{n} (v_\alpha(y_{\alpha,j}(t))-(j-1))
\]
is well-defined (that is, independent of the choice of basis), and each term is $\geq 0$ as $v_\alpha(y_{\alpha,j}(t))\geq j-1$ from the above sequence of inequalities.

If we re-order the $x_j$ so that $v_\alpha(x_1(t)) \leq \dots \leq v_\alpha(x_{n}(t))$ then the 
linear algebra steps  used above to construct the basis  $y_{\alpha,1}(t),\dots,y_{\alpha,n}(t)$  imply that 
\[
 \sum_{j=1}^{n} v_\alpha(y_{\alpha,j}(t)) \geq  \sum_{j=1}^{n} v_\alpha(x_j(t)) ,
 \]
 and so 
 \[
 v_\alpha(T)= \sum_{j=1}^{n} v_\alpha(x_j(t)) -\binom n2 + \Delta_\alpha(x_1,\dots,x_n) \geq   \sum_{j=1}^{n} v_\alpha(x_j(t)) -\binom n2.
\]
One can determine  various lower bounds on 
$\Delta_\alpha(x_1,\dots,x_n)$  (see section \ref{sec: OPP}); for example,  if the minimum of the $ v_\alpha(x_j(t))$ is $0$ and there are $m_\alpha$ values of $j$ with $ v_\alpha(x_j(t))=0$ then  $ \Delta_\alpha(x_1,\dots,x_n)\geq \binom{m_\alpha}2$.

We now study the sequence of degrees of the polynomials $x_j(t)$.  Suppose each  
\[
x_j(t)=a_j t^{e_j}+\cdots \text{ (lower order terms) }
\]
so that  $e_j=\deg x_j$.  If $\deg x_i=\deg x_j$ then
$\deg(a_jx_i-a_ix_j)<\deg x_i$  so we replace the pair $x_i,x_j$ by $x_i,a_jx_i-a_ix_j$. We keep doing this until we end up with a basis 
$y_{\infty,1}(t),\dots,y_{\infty,n}(t)$ for the subspace of $\mathbb C[t]$, generated  by  $x_1(t),\dots,x_{n}(t)$ over $\mathbb C$ where 
\[
0\leq \deg y_{\infty,1}(t)<\deg y_{\infty,2}(t)<\cdots <\deg y_{\infty,n}(t).
\]
This gives the only possible degrees of non-zero elements of $T$ and so   this sequence of integers is unique. Therefore
\[
\deg( T) := \sum_{j=1}^{n} (\deg y_{\infty,j}(t)-(j-1))
\]
is well-defined, and again note that $\deg y_{\infty,j}(t)-(j-1)\geq 0$ for each $j$. Moreover the linear algebra step  implies that 
\begin{equation} \label{eq: Bds on v(x_i),II}
 \sum_{j=1}^{n} \deg y_{\infty,j}(t)\leq  \sum_{j=1}^{n}  \deg x_j(t)  \text{ so that } \deg(T) \leq \sum_{j=1}^{n}  \deg x_j(t)  -\binom n2.
\end{equation}
 We observe a symmetry in the way we make deductions for $\deg$, with all the inequalities reversed from those used or $v_\alpha$.
 Moreover we have the following identity.

\begin{theorem} \label{thm: sumformula for C[t]} Given any  finite dimensional subspace $T$ of $\mathbb C[t]$ over $\mathbb C$
we have 
\begin{equation} \label{eq: Equality for Polyns}
  \sum_{\alpha\in \mathbb C} v_\alpha( T) = \deg(T) .
\end{equation}
\end{theorem}

 \begin{proof} Suppose that $T$ has basis $x_1(t),\dots,x_n(t)$ so that   the Wronskian
\[
W(x_1,  \dots, x_{n}):=\det  \begin{pmatrix} x_1(t) &x_2(t) & \dots & x_{n}(t)     \\  
x_1'(t) &x_2'(t) & \dots & x_{n}'(t)\\
 \dots & \dots & \dots &  \dots\\
x_1^{(n-1)}(t) &x_2^{(n-1)}(t) & \dots & x_{n}^{(n-1)}(t)
 \end{pmatrix}
\]
is non-zero.
We can replace any column by a $\mathbb C$-linear combination of columns (which includes the original column) which changes $W(x_1,  \dots, x_{n})$ by multiplication by a non-zero constant. Therefore the operations which take us from $x_1(t),\dots,x_{n}(t)$ to
$y_{\alpha,1}(t),\dots,y_{\alpha,n}(t)$ in each case above only alter $W(\cdot)$ by a multiplicative constant.
Therefore $\deg(W(x_1,  \dots, x_{n})) = \deg(W(y_{\infty,1},  \dots, y_{\infty,n}))$ and 
\[
v_\alpha(W(x_1,  \dots, x_{n})) = v_\alpha(W(y_{\alpha,1},  \dots, y_{\alpha,n}))  \text{ for each } \alpha\in \mathbb C.
\]

Now suppose that $y_{\infty,j}(t)=c_j t^{d_j}+\cdots$ (lower degree terms) where each $c_j\ne 0$ for $1\leq j\leq n$ and $d_1<\cdots<d_n$.
We claim that the leading terms from the Wronskian will not cancel when we take the determinant, so that 
$\deg(W(y_{\infty,0},  \dots, y_{\infty,n-1})) = \deg(T)$. To see this note that the  $(i,j)$th entry is 
$c_j \cdot \frac{ d_j!}{(d_j-i)!}t^{d_j-i+1}+\cdots$ (lower degree terms). Then   the  leading term of the
 determinant is $\prod_j c_j t^{\sum_j d_j - \binom n2}$ times the determinant of the $n$-by-$n$ matrix with $(i,j)$th entry $\frac{ d_j!}{(d_j-i)!}$. This entry is a (fixed) monic polynomial of degree $i-1$ in $d_j$, so subtracting certain multiples of the previous rows gives a new matrix with the same determinant and with $(i,j)$th entry $d_j^{i-1}$. This is a Vandermonde matrix and has determinant $\pm \prod_{i<j} (d_j-d_i)$ which is non-zero as the $d_i$ are distinct.
 
 An analogous argument shows that  $v_\alpha(W(y_{\alpha,1},  \dots, y_{\alpha,n}))= v_\alpha( T) $
for all $\alpha\in \mathbb C$. We therefore have
\[
W(x_1,  \dots, x_{n}) = c \prod_{\alpha\in \mathbb C} (t-\alpha)^{v_\alpha( T)}
\]
for some non-zero constant $c$. Calculating the degree through this formula yields that
\[
 \sum_{\alpha\in \mathbb C} v_\alpha( T) = \deg(W(x_1,  \dots, x_{n})) =\deg(T). \qedhere
\]
\end{proof}

\subsection{Rational functions and defining the infinite valuation}

\begin{corollary} \label{corr: sumformula for C(t)} Given any  finite dimensional subspace $T$ of $\mathbb C(t)$ over $\mathbb C$
we have 
\begin{equation} \label{eq: Equality for Polyns}
  \sum_{\alpha\in \mathbb C} v_\alpha( T) = \deg(T) .
\end{equation}
\end{corollary} 

\begin{proof}
Given a basis $z_1(t),\dots,z_n(t)$ for $T$, let $h(t)\in \mathbb C[t]$ be the least common denominator of the rational functions
$z_1(t),\dots,z_n(t)$ and let $x_j(t)=h(t)z_j(t)\in \mathbb C[t]$ for each $j$.  Let $T'$ be the space generated by the $x_i(t)$, so that 
$T'=h(t)\cdot T$.  By definition, $v_\alpha(T') = nv_\alpha(h(t))+ v_\alpha(T)$, and $\deg T' = n\deg h + \deg T$.  So the result follows from Theorem \ref{thm: sumformula for C[t]} applied to $T$ and then the space generated by $h(t)$.
\end{proof}

Given the symmetry between the $v_\alpha$ and $\deg$, it makes sense to define 
\[
v_\infty(x(t)):=-\deg x(t).
\]
To unify definitions and inequalities for $v_\alpha$ and $v_\infty$ we now define
\[
v_\infty(T) = -\deg T - 2\binom n2,
\]
so that Corollary \ref{corr: sumformula for C(t)} becomes
\begin{equation} \label{eq: Equality for PolynVals}
  \sum_{\alpha\in \mathbb C\cup \{ \infty\} } v_\alpha( T) =- 2\binom n2 .
\end{equation}

The $n=1$ case is  the usual sum formula \eqref{sec: Extending valuations},  that   the degree of a polynomial is given by the number of its roots counted with multiplicity. So the first interesting case is $n=2$ which we explore in section \ref{sec: OPP}.


\section{General function fields}

The proof in general is the same as  the proof for $\mathbb C[t]$ and $\mathbb C(t)$ with  appropriate modifications.

\subsection{The basic set up}

We will prove a sum formula for finite dimensional subspaces $T$ of $K$, thought of as a vector space over $k$.
First we must define $v(T)$.

 \begin{lemma}  \label{ChangeOfBasis2} 
 Let $T$ be an $n$-dimensional subspace of $K$,    and 
 fix a valuation $v$ of $K$. There exists a basis $y_{v,1},  \dots , y_{v,n}$ of $T$ for which
 \[
0\leq v_\alpha(y_{v,1})<  \cdots < v_\alpha(y_{v,n})  
\]
 We call this an \emph{ordered basis for $T$ with respect to $v$}.
 Moreover for any other ordered basis $y_{v,1}',  \dots , y_{v,n}'$ for $T$ with respect to $v$  we have $v(y_{v,j}')=v(y_{v,j})$ for each $j$.
   \end{lemma}
 
 \begin{proof}  Let $x_1,\dots,x_n$ be a basis for $T$, and express $x_1,\dots,x_n$ as Laurent series in the local parameter $z_v$.
If $v(x_i)=v(x_j)=e$, say then we can write $x_i= c_i z_v^e+\cdots$ and  $x_j= c_j z_v^e+\cdots$  (where ``$+\cdots$'' represents higher order terms, and $c_i,c_j\ne 0$) and we can replace $x_j$ by $X_j:=c_ix_j-c_jx_i$ so that $v(X_j)>e$.
  Iterating until no two terms have the same valuation, we obtain a set $y_{v,1},  \dots , y_{v,n}$ of linearly independent elements of $T$ and so these must be a basis as claimed (as the coefficients lie in a field).
 
 Now if $x\in T$ so that $x=\sum_j c_j y_{v,j}$ then $v(x)=v(y_{v,j_x})$ where $j_x$ is the minimal $j$-value with $c_j\ne 0$.
 Therefore $v(y_{v,1}')\geq v(y_{v,1})$ and, analogously, $v(y_{v,1})\geq v(y_{v,1}')$ so that $v(y_{v,1}')=v(y_{v,1})$.
 Assume that $v(y_{v,j}')=v(y_{v,j})$ for each $j<J$.  Therefore $v(y_{v,J}')>v(y_{v,j}')=v(y_{v,j})$ for all $j<J$, and therefore 
 $v(y_{v,J}')\geq v(y_{v,J})$. Analogously, $v(y_{v,J})\geq v(y_{v,J}')$ so that $v(y_{v,J})=v(y_{v,J}')$.
 \end{proof}

Following Lemma \ref{ChangeOfBasis2} we can define
\[
v(T) := \sum_{j=1}^n v(y_j) - \binom{n}{2}.
\]
The second part of Lemma  \ref{ChangeOfBasis2} shows that this does not depend on the choice of ordered basis.

Suppose that $T$ is a finite dimensional subspace of the ring of integral elements of $K$ (as above, given any finite dimensional subspace $T$ of $K$ we can find an integral elements of $K$ such that $g\cdot T$ is a subset of the ring of integral elements of $K$).

If $v$ is finite then each $v(y)\geq 0$ and so   $v(y_j)\geq j-1$ for each $j\geq 1$. This implies that $v(T)\geq 0$.
If $v$ is infinite then each $v(y)\leq  0$ and so each $v(y_j)\leq j-n$, and therefore $v(T)\leq -2 \binom n2$. We can re-write
\eqref{sum.formula 4 S}  as 
\[
 \sum_{v \text{ finite}}  v(T) = \sum_{v \text{ infinite}}  -v(T) + \binom{n}2 (2g-2).
\]
There are only finitely many infinite places $v$, and so the right-hand side is bounded.  On the other hand, $v(T)$ is a non-negative integer for each finite place and so \medskip

\centerline{$v(T)\ne 0$  for only finitely many valuations $v$.}

\subsection{Local and global derivatives}
If $f\in K$ then $\frac{df}{dz_v}$ is its local derivative at $v$, obtained by formally differentiating the Laurent series of $f$ with respect to $z_v$ so that
$v( \frac{df}{dz_v} ) = v(f)-1$ if $v(f)\ne 0$ and $v( \frac{df}{dz_v} )\geq 0$ if $ v(f)\geq 0$.
Differentiation with respect to $z$ extends uniquely from $k(z)$ to $K$, and is related to the local derivatives by the chain rule,
\[
\frac{df}{dz_v} =  \frac{df}{dz} \cdot \frac{dz}{dz_v}.
\]
Applying \eqref{sum.formula} to $ \frac{df}{dz}\in K$ we find that 
\[
\sum_vv\bigg( \frac{df}{dz_v}\bigg) =\sum_v v\bigg( \frac{df}{dz}\bigg) +\sum_v v\bigg( \frac{dz}{dz_v}\bigg)=\sum_v v\bigg( \frac{dz}{dz_v}\bigg)
\]
 which is independent of the choice of $f\in K\setminus k$; and in fact equals $2g-2$ where $g$ is the \emph{genus} of $K/k$ (see, e.g. \cite[(8)]{Za}). In particular this implies that 
 \begin{equation} \label{sum.for genus}
 \sum_v v\bigg( \frac{dz_v}{dz}\bigg) = - \sum_vv\bigg( \frac{dz}{dz_v}\bigg) = -(2g-2)
\end{equation}

\subsection{Change of variable in the Wronskian}  The Wronskian $W_z(x_1,\dots,x_n)$  is given by
 \[
W_z(x_1,  \dots, x_{n}):=\det  \begin{pmatrix} x_1  &x_2  & \dots & x_{n}     \\  
\frac{dx_1}{dz} &\frac{dx_2}{dz} & \dots & \frac{dx_n}{dz}\\
 \dots & \dots & \dots &  \dots\\
\frac{d^{n-1}x_1}{dz^{n-1}}&\dots & \dots & \frac{d^{n-1}x_n}{dz^{n-1}}
 \end{pmatrix}
\]

\begin{lemma} [Change of variable formula] \label{lem:Change} For any given valuation $v$ we have
\[
 W_{z_v}=W_z \bigg( \frac{dz}{dz_v}\bigg)^{\binom{n}2}.
\]
\end{lemma}

\begin{proof} The chain rule implies that $\frac{df}{dz_v} =  \frac{df}{dz} \cdot \frac{dz}{dz_v}$.
We will show below, through consecutive applications of the chain rule, that for each $i\geq 1$,
 \begin{equation} \label{higher.derivatives}
\frac{d^if}{dz_v^i} =  \sum_{k=1}^i  P_{i,k}(z,z_v) \cdot \frac{d^kf}{dz^k}
\end{equation}
for some functions $P_{i,k}(z,z_v)$ of $z$ and $z_v$ where $P_{i,i}(z,z_v)= (\frac{dz}{dz_v})^i$.  We insert \eqref{higher.derivatives} for $i=0,\dots,n-1$ and $f=x_1,\dots,x_n$ into $W_{z_v}(x_1,  \dots, x_{n})$.
We can subtract a linear combination of previous rows from the $i$th row, so that the $(i,j)$th entry of the resulting matrix is 
\[
\bigg(\frac{dz}{dz_v}\bigg)^i \cdot  \frac{d^if}{dz^i}.
\]
 Dividing $(\frac{dz}{dz_v})^i $ out of the $i$th row for each $i$, we deduce the result claimed.

We prove \eqref{higher.derivatives} by induction on $i\geq 1$. For $i=1$ this is the chain rule, so assume the formula for $i-1$ with $i\geq 2$. Then
\[
\frac{d^if}{dz_v^i} = \frac{d}{dz_v} \bigg( \sum_{k=1}^{i-1} \frac{d^kf}{dz^k} P_{i-1,k}(z,z_v) \bigg)
=  \sum_{k=1}^{i-1}   \frac{d}{dz_v} \bigg( \frac{d^kf}{dz^k} \bigg) \cdot P_{i-1,k}(z,z_v) +  \sum_{k=1}^{i-1} \frac{d^kf}{dz^k}  \frac{dP_{i-1,k}(z,z_v) }{dz_v}  
\]
by the product formula. The chain rule then gives
\[
 \frac{d}{dz_v} \bigg( \frac{d^kf}{dz^k} \bigg)=  \frac{dz}{dz_v} \cdot \frac{d^{k+1}f}{dz^{k+1}}  
\]
and so \eqref{higher.derivatives} follows with 
\[
P_{i,k}(z,z_v) = P_{i-1,k-1}(z,z_v) \cdot \frac{dz}{dz_v}  + \frac{dP_{i-1,k}(z,z_v) }{dz_v}  . \qedhere
\]
\end{proof}

We now apply the sum formula to the Wronskian to obtain the following:

\begin{corollary} [Sum formula for the Wronskian] \label{cor: NewSumFormula}  In a function field of genus $g$ we have
\[
 \sum_v v(W_{z_v}) = \binom{n}2 (2g-2)
\]
\end{corollary}
\begin{proof} 
Now \eqref{sum.formula}, then Lemma \ref{lem:Change} and then \eqref{sum.for genus} give
\begin{align*}
0&=\sum_v  v(W_z)= \sum_v v(W_{z_v}) + \binom{n}2 \sum_v v\bigg( \frac{dz_v}{dz}\bigg)\\
&= \sum_v v(W_{z_v}) - \binom{n}2 (2g-2). \qedhere 
\end{align*}
 \end{proof}

\subsection{Local evaluation of the Wronskian}
 
 \begin{lemma}  \label{ChangeOfBasis} 
  Let $T$ be an $n$-dimensional subspace of $K$ with a basis $x_1,  \dots,  x_{n}$ over $k$.     
 For any  valuation $v$ of $K$ and  ordered basis $y_{v,1},  \dots , y_{v,n}$  with respect to $v$ for $T$, we have
  \[
  v(W_{z_v}(y_{v,1},  \dots , y_{v,n}))=v(W_{z_v}(x_1,  \dots,  x_{n})) .
  \]
 \end{lemma}
 
 \begin{proof}  As before we express $x_1,\dots,x_n$ as Laurent series in the local parameter $z_v$.
If $v(x_i)=v(x_j)=e$, say then we can write $x_i= c_i z_v^e+\cdots$ and  $x_j= c_j z_v^e+\cdots$  (where ``$+\cdots$'' represents higher order terms, and $c_i,c_j\ne 0$) and we can replace $x_j$ by $X_j:=c_ix_j-c_jx_i$ so that $v(X_j)>e$.
This is a linear change of variable and only alters the Wronskian by a constant:
\[
W_{z_v}(x_1,  \dots,X_j ,\dots, x_{n}) = c_i W_{z_v}(x_1,  \dots,x_j ,\dots, x_{n}) 
\]
so that $v(W_{z_v}(x_1,  \dots,X_j ,\dots, x_{n}))=v(W_{z_v}(x_1,  \dots,x_j ,\dots, x_{n}))$. Iterating until no two terms have the same valuation, we obtain a set $y_{v,1},  \dots , y_{v,n}$ of elements of $K$ as claimed.
 \end{proof}

 \begin{lemma}   \label{ValnOfWronskian} 
 Let $v$ be a valuation of $K$ and $y_{v,1},  \dots , y_{v,n}$ be an ordered basis for $T$ with respect to $v$.
  Then
 \[
  v(W_{z_v}(y_{v,1},  \dots , y_{v,n}))= v(T).
 \]
 \end{lemma}
 
  \begin{proof}   
   We write each $y_{v,j}=c_j z_v^{e_j}+\cdots$, with $e_1<e_2<\cdots$ and $c_j\ne 0, c_j\in k$.
   If $v$ is infinite then each $e_j<0$ and so $e_j-(j-1)<0$. If  $v$ is finite  then each $e_j\geq 0$ and so
 each   $e_j-(j-1)\geq 0$ by the hypothesis.
  
 The $(i,j)$th entry of $W_{z_v}(y_{v,1},  \dots , y_{v,n})$ is given by 
 \[
 c_j (e_j(e_j-1)\cdots (e_j-(i-2)) z_v^{e_j-(i-1)}+\cdots
 \]
 Now if we write out the determinant  $W_{z_v}(y_{v,1},  \dots , y_{v,n})$ explicitly as
 \[
 \sum_{\sigma\in S_n} \epsilon(\sigma) \prod_{i=1}^n (W_{z_v}(y_{v,1},  \dots , y_{v,n}))_{i,\sigma(i)}
 \]
 (where $S_n$ is the set of permutations on $n$ letters, and $\epsilon(\sigma)$($=\pm 1$) is the signature of $\sigma$)
 then we see that for each summand
 \begin{align*}
 v\bigg( \prod_{i=1}^n (W_{z_v}(y_{v,1},  \dots , y_{v,n})_{i,\sigma(i)} \bigg) &= \sum_{i=1}^n  v\bigg(  (W_{z_v}(y_{v,1},  \dots , y_{v,n}))_{i,\sigma(i)} \bigg)\\
 & \geq \sum_{i=1}^n  ( e_{\sigma(i)}-(i-1) ) = \sum_{j=1}^n e_j - \binom n2=v(T),
\end{align*}
and so $v(W_{z_v}(y_{v,1},  \dots , y_{v,n}))\geq v(T)$.  To prove this is actually an equality  we need to show that the contributions of the leading terms do not cancel when we sum over all permutations. Going back to our matrix this sum will be given by the determinant of the matrix  formed by the leading terms:
\[
\det  \begin{pmatrix} c_1z_v^{e_1}  & \dots   &c_nz_v^{e_n}       \\  
c_1 e_1 z_v^{e_1-1} &\dots   & c_n e_n z_v^{e_n-1}\\
 \dots & \dots  &  \dots\\
c_1 e_1\cdots (e_1-(n-2)) z_v^{e_1-(n-1)}   & \dots &c_n e_n\cdots (e_n-(n-2)) z_v^{e_n-(n-1)}
 \end{pmatrix}
\]
 \[
= c_1\cdots c_n \cdot z_v^{ \sum_{j=1}^n e_j - \binom n2}\cdot  \det  \begin{pmatrix}1  & \dots   &1     \\  
  e_1   &\dots   & e_n  \\
 \dots & \dots  &  \dots\\
  e_1\cdots (e_1-(n-2))     & \dots &  e_n\cdots (e_n-(n-2))  
 \end{pmatrix} .
\]
 The $(i,j)$th entry in this matrix is a monic polynomial of degree $i-1$ in $e_j$, and so by subtracting multiples of the previous rows, this determinant equals the determinant of the matrix with  $(i,j)$th entry $e_j^{i-1}$. This is a Vandermonde matrix of determinant $\pm \prod_{1\leq i<j\leq n} (e_j-e_i)$ which is $\ne 0$ as all the $e_i$ are distinct. The result follows.
  \end{proof}
 
 \begin{proof}   [Proof of Theorem  \ref{thm: main}]
 Let $x_1,  \dots,  x_{n}$ be a basis for $T$   over $k$.     By Lemma \ref{ChangeOfBasis} we have
 $v(W_{z_v}(x_1,  \dots,  x_{n})) =  v(W_{z_v}(y_{v,1},  \dots , y_{v,n}))$
for any ordered basis $y_{v,1},  \dots , y_{v,n}$ with respect to $v$ for $T$, and this equals $v(T)$
by Lemma \ref{ValnOfWronskian}. The result then follows from Corollary \ref{cor: NewSumFormula}.
    \end{proof}

\section{Understanding $\mathcal E(d_0,\dots,d_n)$}

In any given solution to \eqref{sum x's} and valuation $v$ we may re-organize
the $x_i$ so that if $d_i=v(x_i)$ then $d_0\leq \cdots \leq d_n$, and therefore
\[
\mathcal E(d_0,\dots,d_n)=  \sum_{j=1}^{n} \bigg( d_j - \max_{1\leq i\leq j} (d_{i}-(i-1) ) \bigg) 
\]
We deduce that $\mathcal E(d_0,\dots,d_n)=0$ if and only if $d_1=\dots =d_n$.

We must have 
$d_0=d_1$ in a solution to  \eqref{sum x's} so we can define $\mathcal E'(d_1,\dots,d_n):=\mathcal E(d_0,\dots,d_n)$ without any loss of information.
Moreover  $\mathcal E'(s+d_1,\dots,s+d_n)=\mathcal E'(d_1,\dots,d_n)$ so we select $s=-d_1$ to ensure the first co-ordinate is $0$. Therefore we evaluate $\mathcal E'(0,d_2,\dots,d_n)$.
  
\subsection{Upper bounds on $\mathcal E(d_0,\dots,d_n)$} 
If $d_1=\dots =d_m$ then  $d_j  - \max_{1\leq i\leq j} (d_i-(i-1) ) =0$ for $1\leq j\leq m$; and this is $\leq j-1$ for $m<j\leq n$. Therefore
$\mathcal E'(d_1,\dots,d_n)\leq \binom n2 - \binom m2$.
This yields \cite[Lemma 2]{BM} and \cite[Lemma 1]{Za}. 
Moreover Theorem \ref{corr: sum of the x_i=0} then implies that if  $v(x_{v,m_v})=0$ for each $v\in S$ then
\begin{equation} \label{eq: Zannier.1c}
 H(X) \leq   
  \binom{n}2 (2g-2)+\sum_{v\in S} \bigg( \binom{n}2-\binom {m_v}2\bigg)
\end{equation}
which is a result of Brownawell--Masser \cite[Theorem A]{BM}. Since each $j-1\leq n-1$ we have 
\[
 \binom{n}2-\binom {m_v}2\leq (n-1) \cdot  \#\{ j: v(x_j)>v(x_{v,1})\},
\]
 so that
 \[
 \sum_{v \in S} \bigg( \binom{n}2-\binom {m_v}2\bigg) \leq (n-1)   \sum_{v \in S }\sum_{j:\ v(x_j)>v(x_{v,1})} 1 =(n-1) \sum_{j=0}^n \#\{ v \text{ finite}: v(x_j)>v(x_{v,1})\},
 \] 
and therefore \eqref{eq: Zannier.1c} implies  \cite[Corollary II]{BM}, namely
  \begin{equation} \label{eq: Zannier.7}
H(X)\leq   (n-1)\bigg( (g-1)n+ \sum_{j=1}^n \#\{ v\in S: v(x_j)>v(x_{v,1})\} \bigg) .
\end{equation}
This implies \eqref{eq: Coprime}.

\subsection{Evaluating $\mathcal E'(0,d_2,\dots,d_n)$}\label{sec: E eval}
For $n=2$ we have
 \[
\mathcal E'(0,d_2)  =   \begin{cases} 0 &\text{ if } d_2=0;\\ 1 &\text{ if } d_2\geq 1.\end{cases}
\]
For $n=3$ we have
\[
\mathcal E'(0,d_2,d_3) 
=   \begin{cases}   0 &\text{ if } d_3=0;\\
 1 &\text{ if } d_3=1,d_2=0;\\
2   &\text{ if } d_3=d_2\geq 1;\\
2 &\text{ if } d_3\geq 2, d_2=0;\\
3  &\text{ if } d_3>d_2\geq 1.
\end{cases}  
\]
To get equality in \eqref{sum.formula 4 x} and \eqref{eq: ht ineq2} when $K=\mathbb C[t]$ and $n=3$ we need
\[
D:=\max_i \deg x_i  =   \sum_{v\in S} \mathcal E' (0,a_v,b_v) -6
\]
where the degrees of the $x_i$ are $D\geq D\geq D-a_\infty \geq D-b_\infty$, and the  $v(x_i)$-values are $0\leq 0\leq a_v\leq b_v$.
Can we find any examples that yield equality here?
To simplify this task we could insist that the polynomials are powerful (every root appears to a power at least two), pairwise coprime, and have the same degree $D$, so that we need
$\# \{ \alpha: x_i(\alpha)=0 \text{ for some } i\} =D/2+3$.  For example an identity of the form
\[
r(t)^2x(t)^8 + s(t)^2y(t)^8 + u(t)^2z(t)^8 + v(t)^2w(t)^8=0
\]
where  $x,y,z,w$ of the same degree,   $r,s,u,v$ are linear and $rsuvxyzw$ is squarefree.
 We return to finding examples when $n=3$ in section \ref{sec: uplift of n=2}.

For $n=4$ we have
\[
\mathcal E'(0,d_2,d_3,d_4)  
=   \begin{cases}   
0 &\text{ if } d_4=0;\\
 1 &\text{ if } d_4=1, d_3=0; \\  
  2 &\text{ if } d_4=2, d_3=0; \\  
2 &\text{ if } d_4=  d_3=1,d_2=0;\\    
3 &\text{ if } d_4= 2, d_3=1,d_2=0;\\    
3   &\text{ if } d_4=  d_3=d_2\geq 1;\\  
 3 &\text{ if } d_4\geq 3, d_3=0; \\  
4 &\text{ if } d_4\geq 3, d_3=1,d_2=0;\\    
4   &\text{ if } d_4>d_3=d_2\geq 1;\\  
4 &\text{ if } d_4=d_3\geq 2, d_2=0;\\   
5   &\text{ if } d_4\geq d_3+2, d_3=d_2\geq 1;\\   
5 &\text{ if } d_4>d_3\geq 2, d_2=0;\\  
5  &\text{ if } d_4=d_3>d_2\geq 1;\\  
6  &\text{ if } d_4>d_3>d_2\geq 1. 
\end{cases}  
\]
A simpler description of $\mathcal E'(0,d_2,\dots,d_n)$  might be helpful in applications, but it is far from evident whether there is such a description.

\section{More inequalities}

  \subsection{Using linear combinations} 
  
 Suppose that $x_1,\dots,x_n$ are linearly independent over $k$ and are $S$-units, and let $T$ be the space generated by $x_1,\dots,x_n$ over $k$. Now $S$ must contain all the infinite valuations as well as all the finite valuations for which $\min_i v(x_i)=\min\{ v(t): t\in T\}<0$, so we may assume that if $v\not\in S$ then $v$ is finite and $v(t)\geq 0$ for all $t\in T, t\ne 0$.  Therefore if $u\in T, u\ne 0$ then
\[
 \sum_{v\in S}   \bigg( v(u)-\min_{1\leq i\leq n} v(x_i) \bigg) = -\sum_{v\not \in S}   v(u)-\sum_v \min_{1\leq i\leq n} v(x_i) \leq H(X)\leq \binom{n}2 (2g-2+\# S)
\]
using the sum formula for the sum over $u$, and that $v(x_i)=0$ for all $v\not\in S$, and then by  \eqref{sum.formula 4 H}.
   This is more general than, but inspired by, Zannier's \cite[Theorem 1]{Za}.  
   
Suppose that $x_1,\dots,x_N\in T$ are $S$-units and any $n$ are linearly independent. 
For any $v\in S$ let $v(x_1)\geq v(x_2)\geq \cdots\geq v(x_N)$.  But as $x_1,\dots,x_n$ are linearly independent
we have $v(x_n)=r_v:=\min\{v(t): t\in T\}$ and so $v(x_n)=\cdots=v(x_N)=r_v$. Therefore by \eqref{eq: pure vS2} we have
\[
v(T) +\binom n2 \geq   \sum_{j=1}^{n} v(x_j)  = \sum_{j=1}^{N} v(x_j)  - (N-n) r_v.
\]
 Therefore as  $v(T)\geq 0$ for all $v\not\in S$, \eqref{sum.formula 4 S} implies that
\[
 \binom{n}2 (2g-2+\# S)\geq \sum_{v\in S}  \bigg(   \sum_{j=1}^{N} v(x_j)  - (N-n) r_v \bigg)= (N-n)H(T)
\]
by the sum formula, and as the $x_i$ are all $S$-units, which is  Zannier's \cite[Theorem 2]{Za}.   

\subsection{Applications of \eqref{eq: Zannier.1c} in $ \mathbb C[t]$} \label{sec: Applications}

  In $ \mathbb C[t]$ the only infinite valuation is the degree, in which case $-\min_i v(x_i)=\max_i \deg (x_i)$.
For finite valuations $\min_i v(x_i)=0$ and so $  H(X) = \max_i \deg (x_i)=:D$. Therefore \eqref{eq: Zannier.7} implies
   \[
D\leq   (n-1)\bigg(  \sum_{j=0}^n \#\{ \alpha\in\mathbb C\cup \{ \infty\}: x_j(\alpha)=0\} -n\bigg)
  \]
since there are $\leq n-1$ values of $j$ for which $\deg x_j < D$.

We now show that if $m\geq N^2-1$ then there no   solutions to 
\[
x_0(t)^m+\cdots +x_N(t)^m=0
\]
where the $x_i(t)/x_j(t)\not\in \mathbb C$ whenever $i\ne j$. 
The $x_i(t)^m$ may generate a space of smaller dimension so we reduce to a minimal identity, say
\[
c_0 x_0(t)^m+\cdots +c_nx_n(t)^m=0
\]
with $1<n\leq N$ for non-zero constants $c_j$ and the $x_i(t)$ generate a space of dimension $n$.  We may  assume that 
there is no common factor (else we can divide it out).. If $d=\max_{0\leq i \leq n}\deg x_i$ then $\#\{ \alpha\in\mathbb C\cup \{ \infty\}: x_j(\alpha)=0\} \leq d$ for each $i$ and $D=dm$ so the above gives
   \[
dm<   (n-1)\bigg(  \sum_{j=0}^n \#\{ \alpha\in\mathbb C\cup \{ \infty\}: x_j(\alpha)=0\} \bigg) \leq
(n-1) (n+1) d 
  \]
and so $m<n^2-1\leq N^2-1$. This improves on \cite[Theorem 2.1]{Gy} (which has stronger hypotheses though with the same lower bound on $m$) and improves on \cite[Theorem 2.1]{Gy} which gives $m\geq 16$ when $N=3$.

    It is more complicated to generalize this to 
  \[
c_0z_0(t)^{m_0}+\cdots +c_nz_n(t)^{m_n}=0
\]
since we can have identities like $t+(1-t)=1$ and scale up by powers of $t$ and $(1-t)$. For example,
for any  pairwise coprime integers $a,b,c$ we have 
\[
f(t)^a+g(t)^b=h(t)^c
\]
where
\[
f(t)=t^{\frac{r+1}a}(1-t)^{\frac sa},\ g(t)=t^{\frac{r}b}(1-t)^{\frac {s+1}b}, h (t)=t^{\frac{r}c}(1-t)^{\frac sc}
\]  
and $r$ and $s$ are positive integers for which
$r\equiv -1 \pmod a, 0 \pmod {bc}$ and $s\equiv -1 \pmod b, 0 \pmod {ac}$. We can avoid this type of objection by assuming the $z_j$ have no common factor:

We now show that if
\[
\sum_{i=0}^n\frac  1{m_i}\leq \frac  1{n-1}
\]
then
there is no solution $z_0(t),\dots,z_n(t)\in \mathbb C[t]$ to
\[
c_0z_0(t)^{m_0}+\cdots +c_nz_n(t)^{m_n}=0
\]
where any $n\geq 2$ terms are linearly independent and the $z_i(t)$ have no common factor.
For, if there is a solution take $x_i=c_iz_i^{m_i}$ and observe that
\[
\#\{ \alpha\in\mathbb C\cup \{ \infty\}: x_j(\alpha)=0\} =  \#\{ \alpha\in\mathbb C\cup \{ \infty\}: z_j(\alpha)=0\}\leq 
\begin{cases}
\frac D{m_j} & \text{ if } \deg z_j = \frac D{m_j} ;\\
1+ \frac {D-1}{m_j} & \text{ if } \deg z_j < \frac D{m_j} .
\end{cases}
\]
Now no more than $n-1$ of the $x_j$ have degree $<D$ and so the above equation gives
   \[
D <   (n-1)  \sum_{j=0}^n \frac D{m_j} 
  \]
which implies the result.    This improves on \cite[Theorem 2.3]{Gy}.


 \section{The two polynomials case ($n=2, K=\mathbb C(t)$)} \label{sec: OPP}
  Here we may assume that $T=\langle x(t), y(t)\rangle_{\mathbb C} \subset  \mathbb C[t]$ where $x$ and $y$ have no common factor.  Now if $\alpha\in \mathbb C$ then  $(t-\alpha)$ divides $ax(t)+by(t)\in T$ if and only if 
 $ax(\alpha)+by(\alpha)=0$. Therefore there is a unique such   $(a_\alpha:b_\alpha):=(y(\alpha):-x(\alpha))\in \mathbb P^1(\mathbb C)$. Similarly, the ``root at infinity'' occurs when $\deg (ax(t)+by(t))<D:=\max \{\deg x(t), \deg y(t)\}$, that is the leading coefficient cancels, and this can happen for a unique   $(a_\infty:b_\infty)\in \mathbb P^1(\mathbb C)$.

\subsection{An interpretation of the Riemann-Hurwitz theorem}
 For convenience we will define here $w_\alpha=v_\alpha$ if $\alpha\in \mathbb C$ and $w_\infty=D+v_\infty$ so that $w(z(t))\geq 0$ for all $z\in T$ and all $w$, with
 \[
 \sum_{\alpha\in \mathbb C\cup\{ \infty\} } w_\alpha(ax(t)+by(t))) =D.
 \]
Moreover if $\alpha\in \mathbb C$ then 
 \[
 v_\alpha(T)= v_\alpha(a_\alpha x + b_\alpha y)-1,
 \]
 and 
 \[
  v_\infty(T) = (-\deg(a_\infty x + b_\infty y)-1) + (-D) = -2D+ w_\infty(a_\infty x + b_\infty y)-1.
 \]
 Summing over all $\alpha$ gives
 \begin{align*}
  \sum_v v(T) &= -2D+\sum_{\alpha\in \mathbb C\cup\{ \infty\}} (w_\alpha(a_\alpha x + b_\alpha y)-1)\\
 &  = -2D+\sum_{ a:b \in \mathbb P^1(\mathbb C)} \sum_{\alpha\in \mathbb C\cup\{ \infty\}: x_{a,b}(\alpha)=0} (w_\alpha(x_{a,b})-1)\\
& = -2D+\sum_{ a:b \in \mathbb P^1(\mathbb C)} (D- \# \{\alpha\in \mathbb C\cup\{ \infty\}: x_{a,b}(\alpha)=0\} )
\end{align*}
where $x_{a,b}(t):=ax(t)+by(t)$, and so   \eqref{sum.formula 4 S} implies that 
\[
\sum_{ a:b \in \mathbb P^1(\mathbb C)} (D- \# \{\alpha\in \mathbb C\cup\{ \infty\}: x_{a,b}(\alpha)=0\} = 2D-2.
\]
This can  be obtained directly from the Riemann-Hurwitz formula applied to $\phi(t)=\frac{x(t)}{x(t)+y(t)}$ as in \cite{Si}.
Since $\# \{\alpha\in \mathbb C\cup\{ \infty\}: x_{a,b}(\alpha)=0\} \leq D$ for all $a:b$
one can then deduce that for any finite set $X\subset T$,  
\begin{equation} \label{eq: New Bound}
  \sum_{z(t)\in X} \# \{\alpha\in \mathbb C\cup\{ \infty\}: z(\alpha)=0\}\geq (\# X-2)D+2.
\end{equation}
One can also deduce this from  \eqref{sum.formula 4 x}, though if $\#X>n+1$ then it is more complicated to prove the analogy for $n>2$, the reason being that in the $n=2$ case there is a unique pair $(a:b):=(y(\alpha):-x(\alpha))\in \mathbb P^1(\mathbb C)$ with $ax(\alpha)+by(\alpha)=0$, but this unique pair generalizes to a unique $(n-1)$-dimensional vector space of  $\mathbb P^{n-1}$   for arbitrary $n$.
  
The most interesting case of \eqref{eq: New Bound} is when $\#X=3$ which, by scaling up $x$ and $y$, can be written as
\[
x(t)+y(t)=z(t)
\]
in which case  
\begin{equation} \label{eq: New Bound2}
 \# \{ \text{Distinct roots of } x(t)y(t)z(t) \} \geq  D+2
\end{equation}
where we count $\infty$ as a root if one of $x,y,z$ has degree $<D$.

We can select any M\"obius transformation $t\to \frac{at+b}{ct+d}$ and then define 
\[ X(t):=(ct+d)^D x\bigg(\frac{at+b}{ct+d}\bigg),\ Y(t):=(ct+d)^D y\bigg(\frac{at+b}{ct+d}\bigg),\ Z(t):=(ct+d)^D z\bigg(\frac{at+b}{ct+d}\bigg)\] so that 
$X(t)+Y(t)=Z(t)$.  If we select the transformation so that $-d/c$ is not a root of $x(t)y(t)z(t)$  then $X(t)Y(t)Z(t)$ has degree $3D$ so does not have ``a root at infinity''. This removes the need for the notion of a ``root at infinity''.

\subsection{Examples of equality}

We are interested in understanding the cases of equality in \eqref{eq: New Bound2} where 
\[
 \# \{ \text{Distinct roots of } x(t)y(t)z(t) \} = D+2.
\]
Any such example gives an infinite class of such examples using M\"obius transformations (as in the previous paragraph).
Well-known  examples include  the Pythagorean triple 
\[
f(t)=(2t)^2, g(t)=(t^2-1)^2, h(t)=(t^2+1)^2
\]
 when $D=4$ and,  more interestingly, the Chebyshev polynomials of the first and second kind. These are obtained by taking
\[
f_n(t)+\sqrt{t^2-1} \, g_n(t) = (t+\sqrt{t^2-1})^n
\]
so that 
\[
f_n(t)^2 - (t^2-1) g_n(t)^2 =1
\]
where $\deg f_n=n, \deg g_n=n-1$ and $D=2n$. Note that \eqref{eq: New Bound2} implies that 
\[
2n+2=\deg  ((t^2-1) f_n(t)g_n(t))\geq \# \{ \text{Distinct roots of } (t^2-1) f_n(t)g_n(t)\} \geq  D+2
\]
and so  $f_n$ and $g_n$ must always be squarefree.

My favourite example is 
\[
t(t-2)^3+(1-t)(1+t)^3+(2t-1)^3=0
\]
Here $D=4$ with six roots (including $\infty$). This is interesting because it gives the doubling law for a diagonal cubic:
Suppose that $a,b,c$ are given integers and we have an integer solution $x,y,z$ to $ax^3+by^3=cz^3$. We can obtain a new  integer solution  by taking $t=ax^3/cz^3$ in the above identity, and multipying though by $(cz^3)^4$ to obtain $aX^3+bY^3=cZ^3$ with
\[
X=x(ax^3-2cz^3),\ Y=y(ax^3+cz^3),\ Z=z(cz^3-2ax^3)
\]
A similar example is given by 
\[
t(t^2+6t-3)^4+(1-t)(t^4-28t^3+6t^2-28t+1)^2=(3t^2-6t-1)^4
\]
where $D=9$ and can be used analogously to ``double'' an integral point on $ax^4+by^4=cz^2$.

In the next subsection we describe a methodical way to determine whether examples exist with certain properties, but not 
a guaranteed way how to find them, which is an open question.

\subsection{An equivalent formulation} \label{sec: Belyi}    
If $x(t)+y(t)=z(t)$  let $\phi(t):=x(t)/z(t)$ so  
\[
\phi^{-1}(0)=\{ \alpha: x(\alpha)=0\},\ \phi^{-1}(1)=\{ \alpha: y(\alpha)=0\},\text{ and } \phi^{-1}(\infty)=\{ \alpha: z(\alpha)=0\},\ 
\]
and so $ \# \{ \text{Distinct roots of } x(t)y(t)z(t) \} = D+2$ if and only if the map $\phi:\mathbb P^1 \to \mathbb P^1$ is only ramified over these three points, $\{ 0,1,\infty\}$. This is called a \emph{Bely\u{\i} map} and determining all Bely\u{\i} maps is a very well-known and important problem; see \cite{SV} or \cite{Za}.
Applying methods of algebraic topology to this question gives some practical criteria for finding examples:

Corresponding to each solution is the sequence of multiplicities of zeros which we call the \emph{signature}. Thus
 $(f_1,\dots,f_k); (g_1,\dots,g_\ell); (h_1,\dots,h_m)$ denotes an example of $f(x)+g(x)=h(x)$ where
  $f(x)$ has zeros of multiplicities $f_1,\dots,f_k$ ;
 $g(x)$ has zeros of multiplicities $g_1,\dots,g_\ell$; and
$h(x)$ has zeros of multiplicities $h_1,\dots,h_m$ with 
\[
f_1+\cdots+f_k=g_1+\cdots+g_\ell=h_1+\cdots+h_m=D
\]
and
\[
k+\ell+m=D+2.
\]
 Anything obtained from the examples given by Mobius transformations has this same signature so this is an equivalence class. 
 Applying methods of algebraic topology, in particular the Riemann existence theorem, implies that there is a solution with a given signature
 if and only if there are permutations $\sigma_f, \sigma_g, \sigma_h$ with cycle decompositions of lengths 
 $f_1,\dots,f_k$,  $g_1,\dots,g_\ell$, and $h_1,\dots,h_m$ respectively, such that $\sigma_f \sigma_g \sigma_h$ is the identity, and the subgroup generated by $\sigma_f$ and $ \sigma_g$ is transitive (see \cite{SV} or  \cite[Remark 2, p.96]{Za} for more details). 
 One can show that the lowest degree signature for which there are no examples is
$(1,3);(2,2);(2,2)$.  
 
 Even if these methods imply that there does exist such an identity/Bely\u{\i} map, they do not indicate how to construct an actual example (to quote \cite{SV}, ``Given the description of a Bely\u{\i} map $f$ in the compressed form of a permutation triple,
it has proven difficult in general to determine explicitly an algebraic model for $f$ and the
curve $X$''). We will explore this further in \cite{AI}.

\section{The many polynomials case ($n\geq 3, K=\mathbb C(t)$)} \label{sec: OPP}
 
 We will search for ``best possible'' examples for subspaces of dimension $n\geq 3$ of $\mathbb C[t]$.

\subsection{Best possible examples --- $D$th powers} \label{sec: dth power}

Let $x_j=(t-\alpha_j)^D$ for $1\leq j\leq n$ with $D\geq n-1$ where the $\alpha_j$ are distinct complex numbers, and so
$X=\{ x_1,\dots,x_n \}$ with $T=\langle X \rangle$.  Now   
\[
v_{t-\alpha_j}(X) = D+\Delta( 0,\cdots,0,D) -\binom n2= D+\binom {n-1}2-\binom n2 = D-(n-1)
\]
for $1\leq j\leq n$, and 
\[
\deg(X)=  D+\cdots+D -\Delta( -D,\cdots,-D) -\binom n2=nD - 2\binom n2=n(D-(n-1))
\] 
Therefore
\[
n(D-(n-1))=\sum_{j=1}^n v_{t-\alpha_j}(X)  \leq  \sum_{\alpha\in \mathbb C} v_\alpha( T) = \deg(T) \leq \deg(X)=n(D-(n-1))
\]
which gives equality in \eqref{eq: Equality for Polyns} and therefore in \eqref{sum.formula 4 x}.

One can formulate things more as in \eqref{sum.formula 4 x} by noting that 
\[
v_{\infty}(X) =   -D-\cdots-D +\Delta( -D,\cdots,-D) -\binom n2=-nD
\] 

\subsection{Best possible examples --- Powers and roots of unity} \label{sec:  power + roots}

For $m\geq 2$ and any $k\geq 1$ we have
\begin{align}
\sum_{\zeta: \zeta^m=1} \zeta (x-\zeta)^{km-1} &= \sum_{i=0}^{km-1} (-1)^{i+1}\binom{km-1}i x^{km-1-i} \sum_{\zeta: \zeta^m=1} \zeta^{i+1} \notag \\
&= \sum_{\substack{i=m-1+jm\\ 0\leq j\leq k-1}} (-1)^{i+1}m \binom{km-1}i x^{(k-1-j)m  }
= \sum_{\ell= 0}^{k-1} c_{m,\ell} x^{\ell m} \label{eq: powers2}
\end{align}
for constants $c_{m,\ell}$ taking $\ell=k-1-j$. Here   $n=k+m-1$ and the maximum degree is $km-1$.
Now for 
\[
X=\{ \zeta (x-\zeta)^{km-1} : \zeta^m=1\} \cup \{ x^{\ell m}: 0\leq \ell\leq k-1\}
\]
 we have
\[
v_{x-\zeta}(X)=(km-1)-(n-1)=(m-1)(k-1)
\]
\[
v_{x}(X)=  \sum_{\ell=1}^{k-1}  ( \ell m -(\ell+m-1)) = (m-1)\binom{k-1}2
\]
while 
\[
\deg(X) =\sum_{\ell=0}^{k-1} (\ell m-\ell)   +\sum_{i=1}^{m-1} ((k-1)m+i-(k+i-1)) = (m-1) \binom{k}2+(m-1)^2(k-1) .
\]
 so that 
 \[
 v_{x}(X) + \sum_{\zeta: \zeta^m=1} v_{x-\zeta}(X) = \deg(X) 
 \]
  which gives equality in \eqref{eq: Equality for Polyns} and therefore in \eqref{sum.formula 4 x}. 

\subsection{A not best possible example but still of interest} \label{sec:  power + roots2}
Now
\[
x^m=1+ \prod_{\zeta: \zeta^m=1} (x-\zeta)
\]
so we can rewrite \eqref{eq: powers2} as 
\begin{equation}\label{eq: powers3}
\sum_{\zeta: \zeta^m=1} \zeta (x-\zeta)^{km-1} = \sum_{r=0}^{k-1} C_{m,r} \bigg(   \prod_{\zeta: \zeta^m=1} (x-\zeta) \bigg)^r
\end{equation}
If 
\[
X'=\{ \zeta (x-\zeta)^{km-1} : \zeta^m=1\} \cup \{ C_{m,r}  (   \prod_{\zeta: \zeta^m=1} (x-\zeta))^r: 0\leq r\leq k-1\}
\]
 then
 $v_{x-\zeta}(X')=v_{x-\zeta}(X)$ and $\deg(X')=\deg(X)$ but $v_{x}(X')=0$ so we do not get equality in \eqref{sum.formula 4 x}.

\subsection{Best possible for $n=3$, using best possible for $n=2$} \label{sec: uplift of n=2}

 Taking $m=2$ in \eqref{eq: powers3} gives the identity
\begin{equation}\label{eq: pre abc powers}
(t+1)^{2n-3}-(t-1)^{2n-3} = \sum_{r=0}^{n-2} g_{n,r}((t-1)(t+1))^r.
\end{equation}
Taking $t-1=2a/b$ with $a+b+c=0$ and multiplying through by $-b^{2n-3}$ yields
\begin{equation}\label{eq: abc powers}
a^{2n-3} + b^{2n-3} + c^{2n-3} = \sum_{r=1}^{n-2} h_{n,r}(ac)^r b^{2n-2r-3}
\end{equation}
For now we will take $a,b,c$ to be   polynomials without  common factors, each of degree $D$, such that there are exactly $D+2$ distinct roots of $a(t)b(t)c(t)$. (Later we will take $a,b,c$ to be integers in this construction.)
The right-hand side is divisible by $abc$ but not by $a^2,b^2$ or $c^2$ and so we can write
\begin{equation}\label{eq: abc powers4}
a^k+b^k+c^k=abcd_k
\end{equation}
for any odd integer $k:=2n-3\geq 1$, for some polynomial  $d_k$.  If say $(t-\alpha)^{e_\alpha}\| a(t)b(t)c(t)$, then 
\[
\sum_\alpha (e_\alpha-1)=3D- \sum_\alpha \#\{ \alpha: abc(\alpha)=0\} = 2D-2.
\]
We now consider $X=\{ a^k, b^k, c^k, abcd_k\}\subset T=\langle a^k, b^k, c^k\rangle$.  Therefore if $\alpha\in \mathbb C$ then 
\[
v_\alpha(X) = (ke_\alpha-2)+ (e_\alpha-1)= (k-2)e_\alpha+   3(e_\alpha-1),
\]
so that 
\[
\sum_{\alpha \in \mathbb C} v_\alpha(X) = (k-2)\cdot 3D + 3\cdot (2D-2)=3kD -6.
\]
Now  $v_\infty(X) \geq  -3\cdot (kD)$ and so
\[
-6\leq \sum_v v(X) \leq \sum_v v(T)\leq \binom 32 (2\cdot 0-2) =-6,
\]
and therefore we get equality.  \hfill \qed

\medskip

Here we began with a specific identity   of the form $x_0(t)+x_1(t)+x_2(t)+x_3(t)=0$, and substituted in $t-1=a/b$ where $a,b$ gave a ``best-possible'' example with $n=2$,
to obtain an identity \eqref{eq: abc powers4} which is a ``best-possible'' example with $n=3$.  We might try to generalize this by
starting with an identity of the form $x_0(t)+\cdots+x_n(t)=0$ and substituting in $t-\beta=a/b$ where $a,b$ gives a ``best-possible'' example with $n=2$, and hope for a ``best-possible'' example for arbitrary $n$. One can succeed with the analogous argument provided that for our original $(n+1)$-term identity, there are three  valuations $v$ for which 
\[
v(x_{v,0})= v(x_{v,1})<\dots < v(x_{v,n}).
\]
In \eqref{eq: pre abc powers} this holds for the valuations $v_{-1}, v_{1}$ and $v_\infty$.

\subsection{Infinitely many identities, same subvariety} \label{sec:  Inf Id}
Equation \ref{eq: pre abc powers} implies that 
\[
(t+1)^{2k+1}-(t-1)^{2k+1}=\sum_{j=0}^k c_j(t^2-1)^j ,
\]
a polynomial of degree $k$ in $t^2-1=(t-1)(t+1)$. These terms only have factors $t-1$ and $t+1$.
For $\ell<k/3$ say let $g(t)$ be  \emph{any} polynomial   of degree $\leq 2\ell$ and define
$h(t)=\sum_{j=0}^\ell c_j(t^2-1)^j -g(t)$, 
 to obtain \emph{infinitely many identities}
 \[
 (t+1)^{2k+1}-(t-1)^{2k+1}=\sum_{j=\ell+1}^k c_j(t^2-1)^j +g(t)+h(t)
 \]
  of degree $2k+1$ with $\leq 4\ell+2$ distinct roots.  So here we have infinitely many identities for far fewer roots than the degree; identities that are not related by M\"obius transformations.  However this does not lead to a contradiction to our 
  extended $abc$-conjecture with exponent $1+\epsilon$ ``other than finitely many subvarieties'', because all of these identities can 
 be viewed as  points on the subvariety of $x_0+x_1+\cdots +x_{n}=0$ where $n=k-\ell+3$ with $x_0=-(t+1)^{2k+1}, x_1=(t-1)^{2k+1},$
  \[
x_2^{2k+1}  =c_k^{2k+1}  (-x_0x_1)^{k} \text{ and }  x_{j+1} = -c_kc_{j+1}x_jx_0x_1/c_jx_2^2  \text{ for } j=2,\dots,n-3
  \]
   since $t^2-1=-c_kx_0x_1/x_2$.

 
 \section{Extended abc-conjectures for polynomials and integers}\label{sec:  Extended discussion}

\subsection{Optimal  examples in Browkin and Brzezi\'nski's conjecture} \label{sec:  Optimal BBconj}
Let $a,b,c$ be coprime integers in  \eqref{eq: abc powers}. Then we have a solution to $x_0+\cdots+x_n=0$  in integers where 
\[
\prod_{p|x_0\cdots x_n} p \asymp_k \prod_{p|abc} p.
\]
 Therefore a best possible $abc$-example here yields 
\[
\max_i |x_i| \asymp \bigg( \prod_{p|abc} p\bigg)^{2n-3+o(1)} = \bigg( \prod_{p|x_0\cdots x_n}  p\bigg)^{2n-3+o(1)}
\]
showing that Browkin and Brzezi\'nski's conjecture is best possible if true (our construction is a minor variant of, but more general than, their construction in \cite{BB}).  This, despite the fact that we have shown   this is \emph{not} a ``best possible example'' in our extended $abc$-conjecture in section \ref{sec:  power + roots2}.

\subsection{Browkin and Brzezi\'nski's conjecture for polynomials} \label{sec: BB}

\begin{lemma} \label{lem: select m,n}
Let $F(x,y)\in \mathbb Z[x,y]$ be   homogenous without repeated factors and of degree $\geq 2$. Then there exist arbitrarily large integers $m,n$ such that $F(m,n)\ne 0$ with
$(m,n)=1$ and 
\[
\prod_{p|F(m,n)} p \ll \max\{ |m|, |n|\}^{\deg F-2}.
\]
\end{lemma}

Assuming the $abc$-conjecture this is best possible as \cite[Theorem 5]{Gra} implies that 
\begin{equation} \label{eq: abc lower bd}
\prod_{p|F(r,s)} p \gg \max\{ |r|, |s|\}^{\deg F-2-o(1)}.
\end{equation}

\begin{proof} We will prove this for $G$ which is the product of two linear factors over $\mathbb Z[x,y]$, and for $G$ irreducible. 
The general result then follows  since $F$ must have such a factor $G$, so writing $F=GH$ we have
$F(m,n)=G(m,n)H(m,n)$ and therefore
\begin{align*}
|F(m,n)|&=|G(m,n)|\cdot |H(m,n)| \\
&\ll \max\{ |m|, |n|\}^{\deg G-2} \cdot \max\{ |m|, |n|\}^{\deg H}=\max\{ |m|, |n|\}^{\deg F-2}
\end{align*}
since $|H(m,n)|\ll_H \max\{ |m|, |n|\}^{\deg H}$.  The only possible obstruction is if 
$F(m,n)=0$, in which case $H(m,n)=0$, but this is impossible for $m,n$ coprime and sufficiently large.

Now suppose $F(x,y)= (ax-b y)(Ax-By) $  where $(a,b)=(A,B)=1$.   Let $d=aB-Ab\ne 0$.
For any fixed prime $q\nmid d$ select $k$ to be the largest integer   with $q^k\leq X$.
Let $M=Bq^k-b$ and $N=Aq^k-a$ so that $a M-b N=dq^k$ and $AM-BN=d$.  Therefore $g=(M,N)$ divides $AM-BN=d$, so let 
$m=M/g, n=N/g$ so that $|m|, |n|\asymp_{q,d} X$ with $(m,n)=1$ and therefore
$\prod_{p|(am-bn)(Am-Bn)} p \leq (d/g)q\asymp_{d,q}1$.

Finally suppose $F$ is irreducible of degree $\geq 2$.
Select an integer $r$ such that $F(r,1)$ has a prime factor $q$ which does not divide $r$ nor the discriminant of $F$.
Then we can lift $r$, so that for each $k$ there exists $r_k\equiv r \pmod q$ such that $F(r_k,1)\equiv 0 \pmod{q^k}$ (and $r_{k+1}\equiv r_k  \pmod{q^k}$).
Now let's suppose $(u,v)$ is the  smallest non-zero vector in the lattice $\Lambda_k:=\{ (m,n): m-r_kn \equiv 0 \pmod {q^k}\}$ 
so that $|u|,|v|\leq q^{k/2}$ by the pigeonhole principle.  Now gcd$(u,v)=q^g$, say, with $0\leq g\leq k/2$, so we write
$u=q^gm, v=q^gn$ so that $|m|,|n|\leq q^{k/2-g}\leq q^{K/2}$ with gcd$(m,n)=1$ and $(m,n)\in \Lambda_{K}$, where $K:=k-g\geq k/2$. Therefore $q^K\gg (|m|+|n|)^2$ and so
\[
\prod_{p|F(m,n)} p \leq \frac{|F(m,n)|}{q^{K-1}} \ll \max\{ |m|, |n|\}^{\deg F-2}. \qedhere
\]
\end{proof}

If   \eqref{eq: abc lower bd} holds then $g=o(k)$ in the last part of this proof. Thus the $abc$-conjecture provides information about the $q$-divisibility of the smallest vector in a lattice.

\begin{prop} \label{prop BB}
Assume Browkin and Brzezi\'nski's conjecture. 
For any solution to  \eqref{sum x's} where  the $x_i\in \mathbb Z[t]$ have no common root, generate an $n$-dimensional subspace $T$ and are each of degree $D$, we have
\begin{equation} \label{eq: BB poly}
D \leq  (2n-3)  (\# \{\alpha\in \mathbb C:\  x_i(\alpha)=0 \ \text{ for some } i\} -2) .
\end{equation}
\end{prop}

\begin{proof} Suppose we have a solution to  \eqref{sum x's} where  the $x_i\in \mathbb Z[t]$ have no common root, generate an $n$-dimensional subspace $T$ and are each of degree $D$ with 
\[
D \geq  (2n-3)  (\# \{\alpha\in \mathbb C:\  x_i(\alpha)=0 \ \text{ for some } i\} -2)   +1.
\]
 Let $x_i(t,u):=u^D x_i(t/u)$ so that $\max_i |x_i(r,s)|\asymp \max\{ |r|, |s|\}^D$,
 and let $F(t,u)$ be the product of the distinct irreducible factors of $x_0(t,u)\cdots x_n(t,u)$ so that 
 \[
D \geq  (2n-3)  (\deg F-2)   +1.
\]
 Now by Lemma \ref{lem: select m,n} we deduce that there exist arbitrarily large integers $m,n$ such that  $(m,n)=1$ and 
\[
\prod_{p|x_0(m,n)\cdots x_n(m,n)} p = \prod_{p|F(m,n)} p \ll \max\{ |m|, |n|\}^{\deg F-2}.
\]
 Therefore
\[
\max_i |x_i(m,n)|\gg \max\{ |m|, |n|\}^D\geq  \max\{ |m|,|n|\}^ {(2n-3)  (\deg F-2)   +1} \gg \bigg( \prod_{p|F(m,n)} p \bigg)^{2n-3+\delta+o(1)}
\]
where $\delta= 1/(\deg F-2) $, contradicting Browkin and Brzezi\'nski's conjecture.
\end{proof}

\subsection*{Thanks}
\emph{Many colleagues have made comments on versions of this paper or responded to specific mathematical questions. Thanks to them all especially
Yuri Bilu, Henri Darmon, David Masser (while on a cycling tour of Italy), Carlo Pagano, Joe Silverman, Sameera Vemulapalli, Akshay Venkatesh, Umberto Zannier and Junren Zheng.
I'd like to thank Richard Mason for introducing me to these questions, as a master's student, back in 1983.
I'd also like to thank the folks from Axiom, particularly Ken Ono and Jujian Zhang for verifying the first proof in the paper!}

 \bibliographystyle{plain}

\appendix
\section{AxiomProver's autonomous Lean verification}\label{sec:AxiomProver}

\centerline{ by \emph{Ken Ono and Jujian Zhang} at Axiom Math}
\bigskip

At Axiom Math, we are developing AxiomProver, an AI system for mathematical research that is based on the method of autoformalization. 
This appendix is independent of the rest of the manuscript. 
The text of this entire paper (including this Appendix) is written without the use of AI. 

{Lean is an interactive theorem prover and functional programming language based on dependent type theory, designed to provide a rigorous computational framework for verifying mathematical proofs \cite{Lean}. In this context, autoformalization refers to the automated translation of  natural-language mathematics into machine-verifiable code.}
Lean files are designed to pass a type-checker whereas the natural language paper is meant to convey ideas to human readers.

{Given problems in natural or formal language, AxiomProver attempts to provide a complete formal proof. When successful, the system outputs the files:
\begin{itemize}
  \item \texttt{problem.lean}, a formalization of the problem statement; and
  \item \texttt{solution.lean}, a formalization of an attempted proof
\end{itemize}
}
As a test case, we gave AxiomProver three tasks related to this paper, offering examples of AI assistance in mathematical research.
The relevant files are all posted in the following repository\footnote{%
The formal proofs provided in this work were developed and verified using Lean and mathlib~\textbf{4.26.0}\cite{Mathlib2020}.
Compatibility with earlier or later versions is not guaranteed due to the
evolving nature of the Lean 4 compiler and its core libraries.%
} \footnote{The input to all three tasks contains a configuration file \texttt{.environment} that contains the single line
  ``\emph{lean-4.26.0}.'' This file is henceforth omitted in the appendix}:
\begin{center}
  \url{https://github.com/AxiomMath/Granville}
\end{center}

\subsection*{Task 1 (Autoformalization).} 
We tested AxiomProver's ability to autoformalize a non-formal\footnote{In the formalization community, such proofs are generally called {\it informal} proofs. The distinction between ``formal'' and ``informal'' is about the format, not rigour. In the appendix, we choose the word ``non-formal'' to avoid the connotation of ``informal''.} proof. The task was to autoformalize Section 2 with the input files:
   \begin{itemize}
    \item \texttt{task.md} defines the task.
    \item \texttt{requirement.md} instructs AxiomProver to assume finite dimensionality throughout.
    \item \texttt{LocalGlobal.nterms6.tex} that contains the non-formal proof of Section 2. 
      \end{itemize}
{The system successfully produced a formal statement and a Lean-verified file\footnote{The files can be found in the directory of \url{https://github.com/AxiomMath/Granville/blob/main/Granville/section2/}} corresponding to the non-formal proof.} 
\smallskip

 {We tested AxiomProver's ability to autonomously  refine and formalize existing content as well.}

\subsection*{Task 2 (Sharpening an inequality).} 
We asked AxiomProver to autonomously sharpen (if possible) certain inequalities given in Section~6. The system discovered improvements through a natural case-by-case analysis that a mathematician could have carried out. Theorem A.1 below is this refinement.
 Assume $x_0,x_1,x_2,x_3\in\mathbb{C}[t]$ satisfy $x_0+x_1+x_2+x_3=0$, share no common root, span a $3$-dimensional $\mathbb{C}$-subspace, and have common degree $D:=\deg x_0 = \deg x_1 = \deg x_2= \deg x_3$. For $\alpha\in\mathbb{C}$, let $m(\alpha)$ be the number of $x_0,x_1,x_2,x_3$ vanishing at $\alpha$; necessarily $m(\alpha)\in\{0,1,2\}$. Put $X:=\{\alpha:m(\alpha)\ge 1\}$ and
\begin{itemize}
  \item $N_{1,1}$: points of $X$ where one of $x_0,x_1,x_2,x_3$ has a simple zero,
  \item $N_{1,\ge 2}$: points where one has a zero of order $\ge 2$,
  \item $N_{2,=}$: points where two vanish to the same order,
  \item $N_{2,<}$: points where two vanish to distinct orders;
\end{itemize}
and $N_1:=N_{1,1}+N_{1,\ge 2}$, $N_2:=N_{2,=}+N_{2,<}$.

\begin{theorem}[Sharpened \(n=3\) inequality]\label{thm:sharpened}
\[
 D\le N_{1,1}+2N_{1,\ge 2}+2N_{2,=}+3N_{2,<}-6.
\]
Equivalently, we have
\[
 D\le 2N_1+3N_2-6-(N_{1,1}+N_{2,=}).
\]
\end{theorem}
{The system successfully produced a formal statement and a Lean-verified file\footnote{The files can be found in the directory of \url{https://github.com/AxiomMath/Granville/tree/main/Granville/sharpening_n3/}} corresponding to the non-formal proof.}

\subsection*{Task 3 (Granville's ``simplest configuration'').} 
In Section 6, Andrew Granville invites readers to investigate  the following non-trivial and plausible configuration in the $n=3$ case,
$$
   r^2 x^8 +  s^2 y^8 +  u^2 z^8 +  v^2 w^8 \;=\; 0,
$$
where   $r,s,u,v,x,y,z,w\in\mathbb{C}[t]$, $r,s,u,v$ are linear and $rsuvxyzw$ is squarefree,
as an interesting example for further study.
In response to this suggestion,
AxiomProver autonomously produced and formalized Theorem~\ref{thm:exactfactor} and Corollary~\ref{cor:perfect6} below
with the following input files.
\begin{itemize}
  \item \texttt{task.md} that defines the task.
  \item \texttt{LocalGlobal.nterms6.pdf}: An earlier version of the current manuscript in the pdf format.
\end{itemize}
With the same notation, we assume the following conditions:
\begin{enumerate}[label=(S\arabic*)]
  \item $r^2 x^8 + s^2 y^8 + u^2 z^8 + v^2 w^8 = 0$
  \item $r''=s''=u''=v''=0$ and $x,y,z,w$ all have the same degree;
  \item the pairs $rx,\ sy,\ uz,\ vw$ are pairwise coprime;
  \item $\dim_\mathbb{C}{\operatorname{span}}_\mathbb{C}\{r^2x^8,\,s^2y^8,\,u^2z^8,\,v^2w^8\}=3$;
  \item equality holds in Theorem~\ref{thm:sharpened}.
\end{enumerate}
Write $W:=W_t(r^2x^8,\,s^2y^8,\,u^2z^8)$ for the $3\times 3$ Wronskian, as a polynomial in~$t$, of any three of the four summands (this is well-defined by the above condition).

\begin{theorem}[Exact Wronskian factorization]\label{thm:exactfactor}
{Under the setup above,} there exists a nonzero constant \(\lambda\in \mathbb{C}^\times\) such that
\[
 W=\lambda (xyzw)^6.
\]
\end{theorem}

\begin{corollary}[Perfect-sixth-power determinant]\label{cor:perfect6}
{Assume the setup above. For any $\ell,p\in\mathbb{C}[t]$ with $\ell''=0$, set $F:=\ell^2 p^8$; then $F'=p^6 P_{\ell,p}$ and $F''=p^6 Q_{\ell,p}$, where}
\begin{align*}
 P_{\ell,p}&:=2\ell\ell' p^2+8\ell^2 p p',\\
 Q_{\ell,p}&:=2(\ell')^2p^2+32\ell\ell' p p'+8\ell^2 p p''+56\ell^2 (p')^2.
\end{align*}
Then there exists a nonzero constant \(\mu\in \mathbb{C}^\times\) such that
\[
 \det\begin{pmatrix}
  r^2x^2 & s^2y^2 & u^2z^2\\
  P_{r,x} & P_{s,y} & P_{u,z}\\
  Q_{r,x} & Q_{s,y} & Q_{u,z}
 \end{pmatrix}
 =\mu\, w^6.
\]
\end{corollary}
{The system successfully produced a formal statement and a Lean-verified file\footnote{The files can be found in the directory of \url{https://github.com/AxiomMath/Granville/tree/main/Granville/wronskian/}} corresponding to the non-formal proof.} 
In all three test cases, AxiomProver is able to produce \texttt{problem.lean} files and \texttt{solution.lean} files representing the formal statements and proofs {without additional hint other than those provided in the input files already}. {The test cases offered by this paper illustrate how AxiomProver can assist mathematical research.}

\end{document}